\documentclass[a4paper,10pt]{article}
\usepackage[margin=0.8in]{geometry}
\usepackage{amssymb, amsmath, amsthm, mathrsfs}
\usepackage{cite}
\usepackage{lineno}
\usepackage{tabls}
\usepackage{multirow}
\usepackage{pifont}
\usepackage{graphicx}

\usepackage{color}
\def\qed{\hfill $\Box$\vspace{0.3cm}}

\newtheorem{theorem}{Theorem}[section]%
\newtheorem{lem}[theorem]{Lemma}%
\newtheorem{defi}[theorem]{Definition}%

\def\f{\noindent}

\begin{document}

\title{Paired 3-disjoint path covers of balanced hypercubes}

\author{Mei-Rong Guo \thanks{The work was included in the MS thesis of the author in [Path Cover and Fault-tolerant Cycle Embedding
Analysis of Some Networks, MS Thesis at Beijing Jiaotong University, China, June 2016].}
\hspace{.15in} Rong-Xia Hao \thanks{Corresponding author. This work was partially supported by the National Natural Science Foundation of China (Nos. 11941057 and 11731002) and the 111 Project of China (B16002).}
\hspace{.15in} Mei-Mei Gu\thanks{Current address: Faculty of Mathematics and Physics, Charles University, Prague 11800, Czech Republic. The author was partially supported by OP RDE project No.
CZ.02.2.69/0.0/0.0/16-027/0008495, International Mobility of
Researchers at Charles University and the China Postdoctoral Science Foundation Grant 2018M631322.
\newline\indent\indent
Email address: {\tt 12121625$@$bjtu.edu.cn} (M.-R. Guo), {\tt rxhao@bjtu.edu.cn} (R.-X. Hao), {\tt mmgu@bjtu.edu.cn} (M.-M. Gu)}\\\\
{\small Department of Mathematics, Beijing Jiaotong University, Beijing, 100044, China}
}
%
%
%
\date{}
\maketitle

\begin{abstract}
\baselineskip=12pt
The balanced hypercube $BH_{n}$, proposed by Wu and Huang, is a variation of the hypercube. The paired 1-disjoint path cover of $BH_{n}$ is the Hamiltonian laceability, which was obtained by Xu et al. in [Appl. Math. Comput. 189 (2007) 1393--1401]. The paired 2-disjoint path cover of $BH_{n}$ was obtained by Cheng et al. in [Appl. Math. and Comput. 242 (2014) 127-142]. In this paper, we obtain the paired 3-disjoint path cover of $BH_{n}$ with $n\geq 3$.
This result improves the above known results about the paired $k$-disjoint path covers of $BH_{n}$ for $k=1,2$.

\end{abstract}
\setcounter{page}{1}
\baselineskip=12.5pt


\section{Introduction}
The hypercube network is one of the most popular interconnection networks. The balanced hypercube, proposed by Wu and Huang \cite{Huang2}, is one of the variations of the hypercube. It is a bipartite graph, node-transitive \cite{Huang2} and edge-transitive \cite{J.-X.Zhou}, but has smaller diameter than hypercubes and supports an efficient reconfiguration without changing the adjacent relationship among tasks \cite{Huang2}. In recent years, the balanced hypercube has attracted much attention in the literature \cite{Dongqin}, \cite{R.-X.Hao}, \cite{M.Xu}, \cite{M.-C.Yang}, \cite{J.-X.Zhou}.

One of the most central issues in the study of interconnection networks is to find parallel paths, which is naturally related to routing among nodes and fault tolerance of the network \cite{L.-H.Hsu},\cite{J.A.M.}. Parallel paths correspond to disjoint paths of the graph. The problem of node-disjoint paths has received much attention because of its numerous applications in high performance communication networks, fault-tolerant routings, and so on.

A graph is denoted by $G=(V,E)$, where $V$ is the node set and $E$ is the edge set of $G$. Two nodes $u$ and $v$ are adjacent if $(u,v)\in E$. A path $P[v_0,v_n]$ denoted by $\langle v_0, v_1, \ldots, v_n \rangle$, is a sequence of distinct nodes, where $v_i$ and $v_{i+1}$ are adjacent $(0\leq i\leq n)$. A path $P[v_0, v_n]$ also can be denoted by $\langle v_0, v_1, \ldots, v_i, P[v_i,v_j], v_j, v_{j+1},\ldots v_n \rangle$ , where $P[v_i,v_j]=\langle v_i, v_{i+1}, \ldots, v_j \rangle$. A cycle $C$ is denoted by $\langle v_0, v_1, P[v_1, v_{n-1}],v_0 \rangle$, where $n\geq 3$. A {\it Hamiltonian cycle} (respectively, {\it Hamiltonian path}) of $G$ is a cycle (respectively, path) that traverses every node of $G$ exactly once. A path (resp. cycle) contains $k$ nodes is a {\it $k$-path} (resp. $k$-cycle). A graph is {\it bipartite} if its node sets can be partitioned into two disjoint subsets and each edge with two end nodes from different subsets. A graph G is said to be {\it Hamiltonian connected} if there exists a Hamiltonian path between any two vertices of $G$. It is easy to see that any bipartite graph with at least three vertices is not Hamiltonian connected. For this reason, Simmons \cite{Simmons} introduced the concept of Hamiltonian laceable for Hamiltonian bipartite graphs. A Hamiltonian bipartite graph is {\it Hamiltonian laceable} if there is a Hamiltonian path between any two vertices in different bipartite sets. Obviously, the Hamilton cycle can be embedded in the Hamiltonian connected graphs.
We say that $k$ paths are {\it disjoint} in a network if any two of them have no common node. For an embedding of linear arrays in a network, a cover implies that every node can participate in pipeline computation \cite{Park1}. A $k$-disjoint path cover ($k$-{\it DPC} for short) of a graph is a set of $k$ (internally) disjoint paths that altogether cover every node of the graph.
The problem of disjoint path covers arose from applications of interconnection networks in which the full utilization of nodes is important.  The disjoint path cover problem finds applications in many areas such as software testing, database design, and code optimization \cite{Asdre}, \cite{Ntafos}. We follow \cite{Bondy} for terminologies and notations not specified in the paper.

Let $S$ be a set of $k$ sources $S=\{s_{1},s_{2},\ldots,s_{k}\}$ and a set of $k$ sinks $T=\{t_{1},t_{2},\ldots,t_{k}\}$ in a graph $G$ such that $S\cap T=\emptyset$. A {\it (many-to-many)} $k$-{\it disjoint path cover} of $G$ relate to $S\cup T$ is a set of $k$ disjoint paths joining sources and sinks. It is called {\it paired} if each source $s_{i}$ is joined to a specified sink $t_{i}$. Otherwise, it is called {\it unpaired}.  A {\it paired} (resp. {\it unpaired}) $k$-{\it disjoint path cover} of a bipartite graph $G$ with two distinct partite set $B$ and $W$ is a paired (resp. unpaired) $k$-disjoint path cover of $G$ relate to $S\cup T$ for any two sets $S=\{s_{1},s_{2},\ldots,s_{k}\} \subseteq B$ and $T=\{t_{1},t_{2},\ldots,t_{k}\}\subseteq W$.
 Indeed, if $k=1$, then a paired $1$-disjoint path cover of a network is just a Hamiltonian path connecting the given two nodes. There are some results about $k$-disjoint path cover for $k=1$. However, for $k\geq 2$, the problem of $k$-disjoint path covers of networks is much more difficult. Research attention is focused on some special classes of well-known networks, especially bipartite graphs.
 As examples, for unpaired many-to-many disjoint path cover, Chen, Liu et al. and Park studied a class of bipartite graphs, hypercubes and hypercube-like interconnection networks in \cite{Xie-BinChen4}, \cite{Liu} and \cite{Park3} respectively. Chen, Jo et al. and Lai et al. discussed paired 2-disjoint path covers of hypercubes, bipartite hypercube-like graphs and matching composition network in
 \cite{Xie-BinChen1}, \cite{S.Jo} and \cite{Lai} respectively.

Since the $n$-dimensional balanced hypercube $BH_n$ is a bipartite graph, in the following, its two distinct partite sets are denoted by $B$ and $W$. (In fact, we regard that all vertices in $B$ are given a black color, all vertices in $W$ are given a white color). The paired 1-disjoint path cover of $BH_{n}$, which is the Hamiltonian laceability, was obtained by Xu et al. in \cite{M.Xu}; the paired 2-disjoint path cover of $BH_{n}$ was obtained by Cheng et al. in \cite{Dongqin}. In this paper, we get the paired 3-disjoint path cover of $BH_n$ with $n\geq 3$ which improves the above known results about the paired $k$-disjoint path covers of $BH_{n}$ for $k=1,2$.

\begin{theorem}\label{3-DPC}
Let $BH_{n}$ be the $n$-dimensional balanced hypercube for $n\geq 3$, let $S=\{s_{1}, s_{2}, s_{3}\}\subseteq B$ and $T=\{t_{1}, t_{2}, t_{3}\}\subseteq W$, then $BH_{n}$ has a paired 3-disjoint path cover relative to $S\cup T$, where $B$ and $W$ are two distinct partite sets of $BH_{n}$.
\end{theorem}

The rest of this paper is organized as follows. Section 2 introduces definitions and preliminaries of the balanced hypercube. The proof of our main result is given in Section 3. Section 4 concludes this paper.

\section{Definitions and preliminaries}
In this section, we give definitions and basic properties of the balanced hypercube.

An $n$-dimensional balanced hypercube \cite{Huang2}, denoted by $BH_{n}$, is defined as follows. In the following, $"+"$ and $"-"$ always mean an operation with mod 4. For convenience, ``mod 4" is omitted.

\begin{figure}[ht]
\begin{center}
\unitlength 4mm
\begin{picture}(20,12)

\put(1, 3){\circle*{0.4}} \put(0.6, 2){$3$}

\put(1,3){\line(1, 0){5}} \put(1,3){\line(0, 1){5}}

\put(6, 3){\circle{0.4}}\put(5.7, 2){$2$}

\put(6,3){\line(0, 1){5}}\put(1,8){\line(1, 0){5}}

\put(1, 8){\circle{0.4}} \put(0.6, 8.3){$0$}

\put(6, 8){\circle*{0.4}}\put(5.7, 8.3){$1$}

\put(10, 1){\circle*{0.4}}\put(10,1){\line(0, 1){9}}\put(10,1){\line(1, 0){3}}\put(10,1){\line(2, 1){6}}
\qbezier(10, 1)(14.5, -1)(19, 1)\put(8, 1){{\footnotesize $(3,3)$}}\put(8, 4){{\footnotesize $(0,3)$}}
\put(8, 7){{\footnotesize $(3,0)$}}\put(8, 10){{\footnotesize $(0,0)$}}

\put(10, 4){\circle{0.4}}\put(10,4){\line(1, 0){9}}\put(10,4){\line(1, 2){3}}

\put(10, 7){\circle*{0.4}}\put(10,7){\line(1, 0){9}}\put(10,7){\line(1, -2){3}}

\put(10, 10){\circle{0.4}}\put(10,10){\line(2, -1){6}}\put(10,10){\line(1, 0){3}}
\qbezier(10, 10)(14.5, 12)(19, 10)

\put(13, 1){\circle{0.4}} \put(13,1){\line(0, 1){3}}\qbezier(13, 1)(14.5, 5.5)(13, 10)
\put(12.5, 0){{\footnotesize $(2,3)$}}
\put(15, 0){{\footnotesize $(3,2)$}}

\put(13, 4){\circle*{0.4}}\put(13,4){\line(2, -1){6}}\put(12.5, 4.2){{\footnotesize $(1,3)$}}
\put(15, 4.2){{\footnotesize $(0,2)$}}

\put(13, 7){\circle{0.4}}\put(13,7){\line(0, 1){3}} \put(13,7){\line(2, 1){6}}\put(12.5, 6){{\footnotesize $(2,0)$}}
\put(15, 6){{\footnotesize $(3,1)$}}

\put(13, 10){\circle*{0.4}} \put(12.5, 10.5){{\footnotesize $(1,0)$}}\put(15, 10.5){{\footnotesize $(0,1)$}}

\put(16, 1){\circle*{0.4}}\put(16,1){\line(0, 1){3}}\qbezier(16, 1)(14.5, 5.5)(16, 10)\put(16,1){\line(1, 2){3}}

\put(16, 4){\circle{0.4}}

\put(16, 7){\circle*{0.4}}\put(16,7){\line(0, 1){3}}

\put(16, 10){\circle{0.4}}\put(16,10){\line(1, 0){3}} \put(16,10){\line(1, -2){3}}

\put(19, 1){\circle{0.4}}\put(19,1){\line(0, 1){9}}\put(19,1){\line(-1, 0){3}}

\put(19, 4){\circle*{0.4}}

\put(19, 7){\circle{0.4}}

\put(19, 10){\circle*{0.4}}\put(19.2, 1){{\footnotesize $(2,2)$}}\put(19.2, 4){{\footnotesize $(1,2)$}}
\put(19.2, 7){{\footnotesize $(2,1)$}}\put(19.2, 10){{\footnotesize $(1,1)$}}

\end{picture}
\end{center}\vspace{-0.3cm}

\caption{{\ $BH_1$ and $BH_2$ }}\label{Fig-1}

\end{figure}
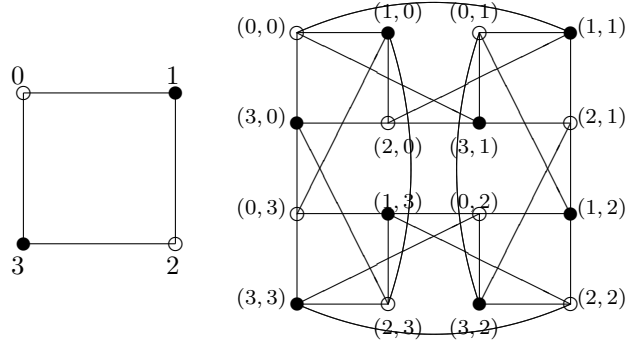

\begin{defi}\label{defi-1}
For $n\geq1$, $BH_n $ consists of $2^{2n}$ nodes $(a_0, a_1, \ldots,a_{i-1}, a_i, a_{i+1}, \ldots, a_{n-1})$, where $a_i\in\{0,1,2,3\}$ $(0\leq i\leq n-1)$. Every node $(a_0,a_1,\ldots,a_{i-1}, a_i, a_{i+1},\ldots, a_{n-1})$ connects the following $2n$ nodes:
\begin{enumerate}
  \item [{\rm (1)}] $(a_0\pm1, a_1, \ldots, a_{i-1}, a_i, a_{i+1}, \ldots, a_{n-1}) $ and
  \item [{\rm (2)}] $(a_0\pm1, a_1, \ldots, a_{i-1}, a_i+(-1)^{a_0}, a_{i+1}, \ldots, a_{n-1})$, where $i$ is an integer with $1\leq i\leq n-1$.
\end{enumerate}
\end{defi}
$BH_1$ and $BH_2$ are illustrated in Figure~\ref{Fig-1}.

In $BH_n$, $a_0$ of a node $(a_0, a_1,\ldots,a_{i-1}, a_i, a_{i+1}, \ldots, a_{n-1})$ is named $inner$ $index$, and $a_i$ for each $i$ $(1\leq i\leq {n-1})$ are named $i$-$dimensional$ $index$. Clearly $BH_n$ is a bipartite graph. In what follows, we refer to a node with an odd inner address as a black node and a node with an even inner address as a white node. The set of all black nodes, say  $B$, and the set of all white nodes, say $W$, make the desired partition.

Based on Definition~\ref{defi-1}, the equivalent definition of $BH_n$ is as follows.

\begin{defi}\label{defi-2}
$BH_n $ is recursively constructed as follows:
\begin{enumerate}
\vspace{-0.2cm}
  \item [{\rm (1)}] $BH_1 $ is a cycle consisting of four nodes labeled as 0,1,2,3 respectively.
 \vspace{-0.2cm}
  \item [{\rm (2)}] For $n\geq2 $, $BH_n$ consists of four copies of $BH_{n-1}$, denoted by $BH_{n-1}^i$, for every integer $i $ with $0\leq i\leq3$. Every node $(a_0, a_1,\ldots, a_{n-2},i)$ of $BH_{n-1}^i$ connects to two extra nodes:\\
(2.1) $(a_0\pm1, a_1,\ldots, a_{n-2},i+1)$ in $BH_{n-1}^{i+1}$, if $a_0$ is even.\\
(2.2) $(a_0\pm1, a_1,\ldots, a_{n-2},i-1)$ in $BH_{n-1}^{i+1}$, if $a_0 $ is odd.
\end{enumerate}
\end{defi}

We classify the edges of $BH_n $ as follows. If two adjacent nodes $u$, $v$ differ in only the inner index, the edge $(u,v)$ is said to be $0$-$dimensional$ and $v$ is a $0$-dimensional neighbor of $u$. Likewise if two adjacent nodes $u$, $v$ not only differ in the inner index, but also differ in some $i$-dimensional index $(1\leq i\leq {n-1})$, the edge $(u,v)$ is said to be $i$-$dimensional$ and $v$ is the {\it $i$-dimensional neighbor} of $u$. Let $E_i$ for $i\in \{0,1,2,\ldots,n-1\}$ denote the set of all edges of $i$-dimensional edges, $E(BH_n)=\bigcup_{i=0}^{n-1}E_i$. For $i\in\{0,1,2,3\}$ and $1\leq j\leq {n-1}$, we use $BH_{n-1}^{j,i}$ to denote $(n-1)$-dimensional sub-balanced hypercubes of $BH_n$ induced by all vertices labeled by $(a_0,a_1,\ldots,a_{j-1},i,a_{j+1},\ldots,a_{n-1})$. Obviously, $BH_n-E_j=\bigcup_{i=0}^3 BH_{n-1}^{j,i}$ and $BH_{n-1}^{j,i}\cong BH_{n-1}$. If $j=n-1$, $BH_{n-1}^{j,i}$ and $E_j$ are denoted by $BH_{n-1}^i$ and $E_c$, respectively, where $i\in\{0,1,2,3\}$, then we have $BH_n-E_c=\bigcup_{i=0}^3 BH_{n-1}^i$. The edge in the set $E_c$ is named the crossing edge.
For any node $u=(a_0,a_1,\ldots, a_{n-1})\in V(BH_n)$, if $v=(a_0+2,a_1,\ldots, a_{n-1})$, we say $u$ and $v$ are {\it a symmetric pair} or $u$ is {\it symmetric with} $v$.

\begin{lem}\label{node-trans} {\rm \cite{Huang2}}
The balanced hypercube $BH_{n}$ is node transitive.
\end{lem}

\begin{lem}\label{H-l} {\rm\cite{M.Xu}}
The balanced hypercube $BH_{n}$ is Hamiltonian laceable.
\end{lem}

\begin{lem}\label{8-cycle} {\rm\cite{M.Xu}}
Let $(u,v)$ be an edge of $BH_{n}$. Then $(u,v)$ is contained in an $8$-cycle in $BH_{n}$ such that $E(C)\cap E(BH_{n-1}^i)=1$, where $i=0,1,2,3$.
\end{lem}

\begin{lem}\label{2-DPC} {\rm\cite{Dongqin}}
Let $S=\{s_{1},s_{2}\}\subseteq B$ and $T=\{t_{1},t_{2}\}\subseteq W$ be two distinct partite sets of $BH_{n}$, then $BH_{n}$ has a paired many-to-many $2$-disjoint path cover relative to $S\cup T$, where $n\geq 1$.
\end{lem}

\begin{lem}\label{4-copies}
Let $t_{1},t_{2},t_{3}$ be three distinct white nodes in $BH_{n}$ for $n\geq2$, then there exists an $l\in\{0,1,\ldots,n-1\}$ such that $BH_{n}$ is divided into four copies of $BH_{n-1}$, say $BH_{n-1}^{l,i}$, for $i=0,1,2,3$ along the $l$-dimension such that $t_{1},t_{2},t_{3}$ belong to at least two $BH_{n-1}^{l,i}$'s.
\end{lem}
\f {\bf Proof.}
For any three white nodes $t_{1},t_{2},t_{3}$ in $BH_{n}$, assume $t_j=(a_0^j,a_1^j,\ldots,a_{n-1}^j)$, $j\in \{1,2,3\}$, there exists $l\in\{0,1,\ldots,n-1\}$ such that at least two of $a_l^1$, $a_l^2$ and $a_l^3$ are not equal, otherwise $t_{1}$, $t_{2}$ and $t_{3}$ are the same node, which is a contradiction. Then $BH_{n}$ is divided into four copies of $BH_{n-1}$, say $BH_{n-1}^{l,i}$'s along $l$ , for $l\in\{1,2\}$, such that $BH_n-E_l=\bigcup_{i=0}^3 BH_{n-1}^{l,i}$ and $t_{1},t_{2},t_{3}$ belong to at least two $BH_{n-1}^{l,i}$'s, $i\in \{0,1,2,3\}$.\hfill\qed

\begin{lem}\label{5-path}
 For any black node $s$ in $BH_{n}$ with $n\geq 2$, there exist two $5$-paths, say $\langle s,a_{t},b_{t},c_{t},d_{t}\rangle$, passing $s$ and
  Hamiltonian cycles $C_t$ in $BH_{n}$ such that $C_t$ contains the path $\langle s,a_{t},b_{t},c_{t},d_{t}\rangle$, where $\langle a_{1},b_{1},c_{1},d_{1}, a_{1}\rangle$ and $\langle a_{2},b_{2},c_{2},d_{2},a_{2}\rangle$ are two node-disjoint 4-cycles, $a_{t}$ and $b_{t}$ are symmetric with $c_{t}$ and $d_{t}$, respectively, for $t=1,2$.
\end{lem}
\f {\bf Proof.}
We prove the lemma by induction on $n$.

For $n=2$, by Lemma~\ref{node-trans}, we may assume that $s=(1,0)$, then there exists a $5$-path $P_1=\langle (1,0),(0,0),(1,1)$, $(2,0),(3,1)\rangle$ and a Hamiltonian cycle
$C_1=\langle (1,0),(0,0),(1,1),(2,0),(3,1),(0,1),(1,2),(2,1),(3,2),(2,2),(1,3)$, $(0,2),
(3,3),(2,3),(3,0),(0,3),(1,0)\rangle$ which contains $P_1$. The other $5$-path is $P_2=\langle (1,0),(0,3),(3,3),(2,3),(1,3)\rangle$ and $C_2=\langle (1,0),(0,3),(3,3),(2,3),(1,3),(2,2),(3,2),(0,2),(1,2),(2,1),
(3,1),(0,1),(1,1),(2,0),(3,0),(0,0),(1,0)\rangle$ is the desired Hamiltonian cycle which contains $P_2$. By the definition of $BH_2$, we can see that $\langle (0,0),(1,1),(2,0),(3,1)$, $(0,0)\rangle$ and $\langle(0,3),(3,3),(2,3),(1,3),(0,3)\rangle$ are $4$-cycles, $(0,0)$, $(1,1)$, $(0,3)$ and $(1,3)$ are symmetric with $(2,0)$, $(3,1)$, $(2,3)$ and $(3,3)$, respectively.

Assume that the lemma is true for $BH_{n-1}$ for $n\geq 3$. We consider $BH_{n}$ as follows.

By the Definition~\ref{defi-2}, $BH_{n}$ can be divided into four components $BH_{n-1}^i$'s along $(n-1)$-dimension with $i=0,1,2,3$. For any black node $s\in V(BH_{n})$, without loss of generality, we may assume that $s\in V(BH_{n-1}^0)$. By the inductive hypothesis, there exist two $5$-paths $\langle s, a_{t},b_{t},c_{t},d_{t}\rangle$ in $BH_{n-1}^{0}$ and two  Hamiltonian cycles, say $C_t^0$, in $BH_{n-1}^0$ such that $C_t^0$ contains the path $\langle s,a_{t},b_{t},c_{t},d_{t}\rangle$, where $a_{t}$ and $b_{t}$ are symmetric with $c_{t}$ and $d_{t}$, respectively, $\langle a_{t},b_{t},c_{t},d_{t}, a_t\rangle$ are two node-disjoint $4$-cycles in $BH_{n-1}^{0}$ for $t=1,2$.
Let $(u_t,v_t)$ be an edge in $C_{t}^0$, $u_t$ be white and $u_t,v_t \notin \{s,a_{t},b_{t},c_{t},d_{t}\}$. Let $P_t=C_t^0-(u_t,v_t)$, then $P_t$ contains the path $\langle s,a_{t},b_{t},c_{t},d_{t}\rangle$. By Lemma~\ref{8-cycle}, $(u_t,v_t)$ is contained in an $8$-cycle $C_t'$ in $BH_{n}$ such that $|E(C_{t}')\cap E(BH_{n-1}^i)|=1$ with $i=0,1,2,3$. Let $E(C_{t}')\cap E(BH_{n-1}^i)=(u_t^i,v_t^i)$, by Lemma~\ref{H-l}, there exists an Hamiltonian path $P_t^i=\langle u_t^i,v_t^i \rangle$ in $BH_{n-1}^i$ for each $i\in \{1,2,3\}$. Then $C_t=\langle v_t,P_t,u_t,v_t^1,P_t^1,u_t^1,v_t^2,P_t^2,u_t^2,v_t^3,P_t^3,u_t^3,v_t\rangle$ which contains the path $\langle s,a_{t},b_{t},c_{t},d_{t}\rangle$ for $t=1,2$ are our desired Hamiltonian cycles. The proof is complete. \hfill\qed

\begin{lem}\label{BH2}
Let $S=\{s_{1},s_{2},s_{3}\}\subseteq B$ and $T=\{t_{1},t_{2}\}\subseteq W$, where $B$ and $W$ are two distinct partite sets of $BH_{2}$, then we can find a node $t_{3}\in W-\{t_{1},t_{2}\}$ such that $BH_{2}$ has a paired $3$-disjoint path cover relative to $S\cup (T\cup\{t_{3}\})$.
\end{lem}
\f {\bf Proof.} Since $BH_2$ is node-transitive, we may assume that $t_1=(0,0)$. Note that the following nodes are symmetric to each other \cite{M.-C.Yang}: $(2,i)$ and $(0,i)$, $(3,i)$ and $(1,i)$ for $i=0,1,2,3$. So we may assume that $t_2\in\{(2,0),(0,1),(0,2),(0,3)\}$,  $s_1\in\{(1,0),(1,1),(1,2),(1,3)\}$ and $s_2\in\{(3,0),(3,1),(3,2),(3,3)\}$. For any given quadruples $t_1,t_2,s_1$ and $s_2$, we consider every black node that is different from $s_1$ and $s_2$. For each case, we can find $t_{3}\in W-\{t_{1},t_{2}\}$ such that $BH_{2}$ has a paired $3$-disjoint path cover relative to $S\cup (T\cup\{t_{3}\})$. All cases are listed in the appendix tables.\hfill\qed

Recall that $B$ and $W$ are two distinct partite sets of $BH_{n}$, $S=\{s_{1},s_{2},s_{3}\}\subseteq B$ and $T=\{t_{1},t_{2},t_{3}\}\subseteq W$. For each pair of $\{s_j,t_j\}$ with $s_j\in V(BH_{n-1}^g)$ and $t_j\in V(BH_{n-1}^h)$ for $j\in \{1,2,3\}$ and $g,h\in \{0,1,2,3\}$, let

$$M_{j}^i=\Bigm\{ \begin{matrix}
            1, &i=g,g+1,\ldots,h,~~~~~~~~~~~\\
            0, &i=h+1,h+2,\ldots,g-1,\\
          \end{matrix} $$ where $i\in \{0,1,2,3\}$, $``+"$ is mod $4$.

The value of $M_{j}^i$ depends on the location of $\{s_j,t_j\}$. For example, if $s_j\in V(BH_{n-1}^2)$, $t_j\in V(BH_{n-1}^0)$, then $M_{j}^2=M_{j}^3=M_{j}^0=1$ and $M_{j}^1=0$.

Furthermore, for each $i\in \{0,1,2,3\}$, let $\beta_{i}=M_{1}^i+M_{2}^i+M_{3}^i$. By the definition of $M_{j}^i$, we have $0\leq\beta_{i}\leq3$, $i\in \{0,1,2,3\}$. We use $f_k$ to represent the number of $\beta_{i}$'s with $\beta_{i}=k$, $k\in \{0,1,2,3\}$. Obviously, $0\leq f_k\leq 4$ for each $k\in \{0,1,2,3\}$ and $f_0+f_1+f_2+f_3=4$.

\section{Proof of the Theorem 1.1}

\f {\bf Proof.}
We prove the theorem by induction on $n$. The proof for the  basis, i.e., when $n=3$, is similar to that for the induction step, i.e., $n\geq 4$. We combine these proofs to avoid repetition. Assume that the theorem is true in $BH_{n-1}$ for $n\geq 4$. We consider $BH_{n}(n\geq 4)$ as follows.

By Lemma~\ref{4-copies}, there exists a partition along a dimension $l$ ($l\in\{0,1,\ldots,n-1\}$) such that $t_{1},t_{2},t_{3}$ are not in the same $BH_{n-1}^{l,i}$ for $i\in\{0,1,2,3\}$. Without loss of generality, we may assume $l=n-1$, i.e., $BH_n-E_{n-1}=\bigcup_{i=0}^3 BH_{n-1}^i$. We consider the following two cases.

{\bf Case 1.} $f_0=0$ or $f_3<2$.

For each pair ${s_j,t_j}$, $s_{j}\in V(BH_{n-1}^g)$ and $t_{j}\in V(BH_{n-1}^h)$ with $j\in \{1,2,3\}$, we choose nodes as follows (recall that the operation is mod 4):
\begin{enumerate}
  \item [{\rm (1)}]When $h-g=1$, we choose a white node, say $a_{j}^g$, in $BH_{n-1}^g$ and a black node, say $b_{j}^h$, in $BH_{n-1}^h$, such that $(a_{j}^g,b_{j}^h)\in E_{c}$.
  \item [{\rm (2)}]When $h-g=2$, we choose two white nodes $a_{j}^g$ and $a_{j}^{g+1}$ in $BH_{n-1}^g$ and $BH_{n-1}^{g+1}$ respectively, two black nodes $b_{j}^{g+1}$ and $b_{j}^h$ in $BH_{n-1}^{g+1}$ and $BH_{n-1}^h$, respectively, such that $(a_{j}^g,b_{j}^{g+1}),(a_{j}^{g+1},b_{j}^h)\in E_{c}$;
  \item [{\rm (3)}]When $h-g=3$, we choose three white nodes $a_{j}^g,a_{j}^{g+1}$ and $a_{j}^{g+2}$ in $BH_{n-1}^g$, $BH_{n-1}^{g+1}$ and $BH_{n-1}^{g+2}$ respectively, three black nodes $b_{j}^{g+1},b_{j}^{g+2}$ and $b_{j}^h$ in $BH_{n-1}^{g+1}$, $BH_{n-1}^{g+2}$ and $BH_{n-1}^h$ respectively, such that $(a_{j}^g,b_{j}^{g+1}),(a_{j}^{g+1},b_{j}^{g+2}),(a_{j}^{g+2},b_{j}^h)\in E_{c}$.
\end{enumerate}

\f Let $A$ (resp. $B$) be the set of white (black) nodes chosen by (1)-(3). Since $n\geq 3$, $\mid V(BH_{n-1}^i)\mid\geq2^3$, we can assume $S\cap T\cap A \cap B=\emptyset$. The nodes in $S\cup T\cup A \cup B$ are paired in $BH_{n-1}^i$. In other words, if there is one black node in $(S\cup T\cup A \cup B)\cap BH_{n-1}^i$, there must be one corresponding white node in the $(S\cup T\cup A \cup B)\cap BH_{n-1}^i$ for $i\in\{0,1,2,3\}$. By the definition of $\beta_{i}$, we can see that every $BH_{n-1}^i$ for each $i\in\{0,1,2,3\}$ has $\beta_{i}$ pairs of nodes in $S\cup T\cup A \cup B$ and each pair of nodes contains only one white node and one black node. We deal with the following cases.

{\bf Subcase 1.1.} $f_0=f_3=0$.

\begin{figure}
\begin{minipage}[t]{0.5\linewidth}
\centering
\includegraphics[height=5cm,width=7cm]{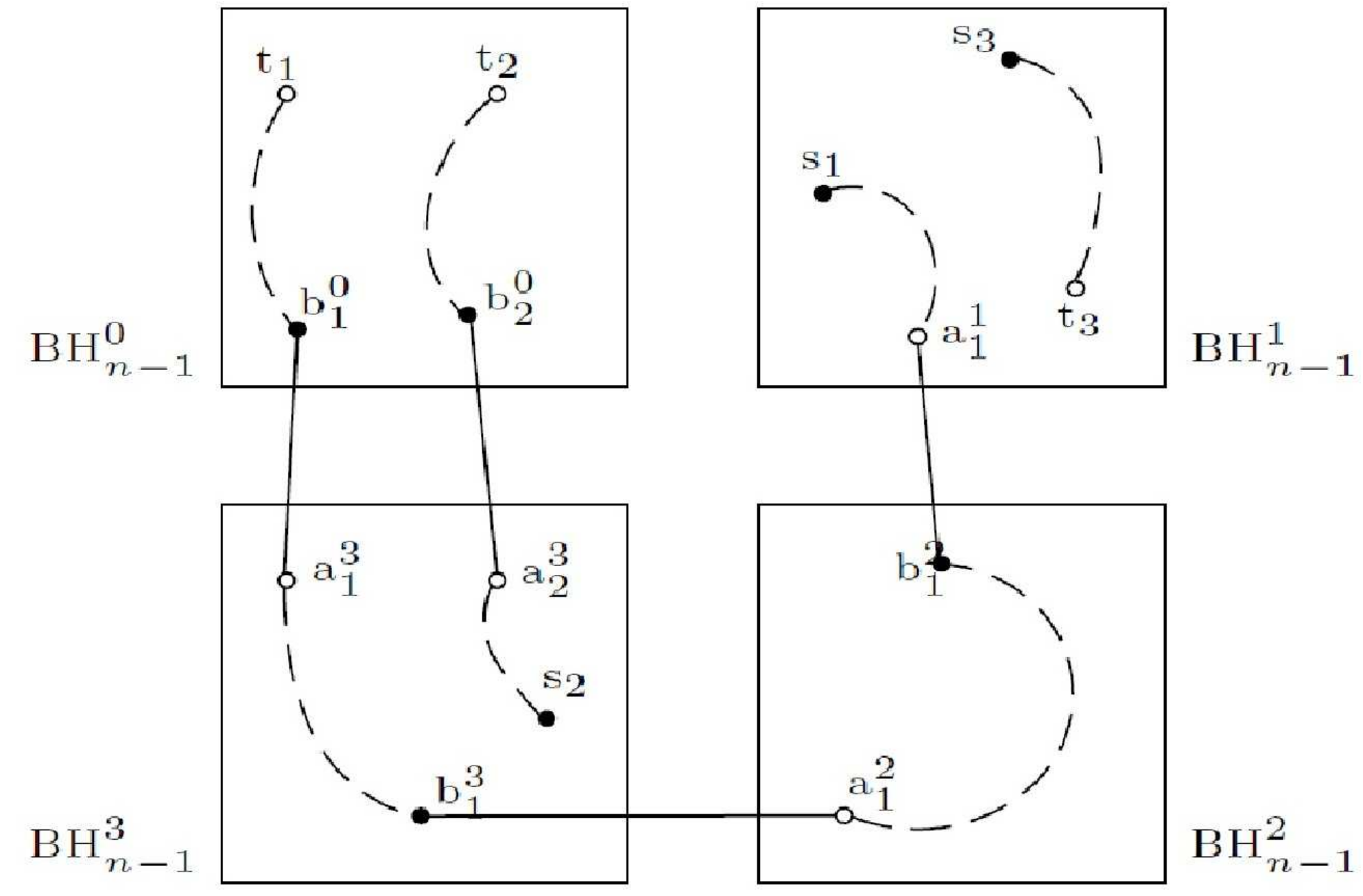}
\vskip-0.3cm
\caption{Illustration of 1.1. in Theorem~\ref{3-DPC}}
\label{Fig-3}
\end{minipage}
\begin{minipage}[t]{0.5\linewidth}
\centering
\includegraphics[height=5cm,width=7cm]{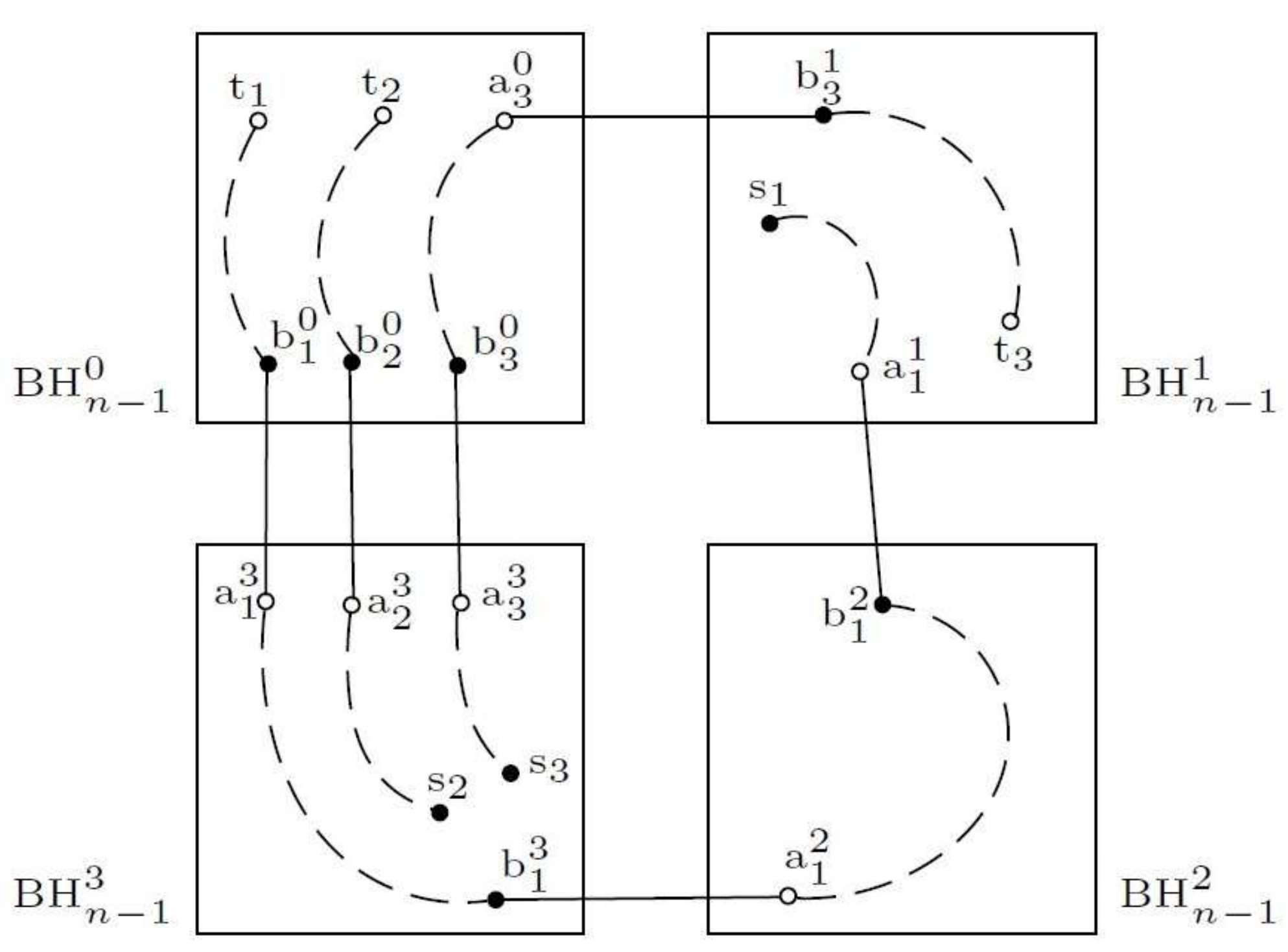}
\vskip-0.3cm
\caption{Illustration of 1.2. in Theorem~\ref{3-DPC}}
\label{Fig-4}
\end{minipage}
\end{figure}

It implies that $f_1+f_2=4$, i.e., for each $i\in \{0,1,2,3\}$, $BH_{n-1}^i$ has $1\leq \beta_{i}\leq 2$ pairs in $S\cup T\cup A \cup B$. Then $(f_1,f_2)\in\{(1,3),(3,1),(0,4),(4,0),(2,2)\}$. Since the discussions are similar, we only consider the situation that $(f_1,f_2)=(1,3)$ and assume $s_{1},s_{3}\in V(BH_{n-1}^1),s_{2}\in V(BH_{n-1}^3)$ and $t_{1},t_{2}\in V(BH_{n-1}^0),t_{3}\in V(BH_{n-1}^1)$ (see Figure~\ref{Fig-3}). Recall that for each $i\in \{0,1,2,3\}$, $BH_{n-1}^i$ has $\beta_{i}$ pairs of black and white nodes in $S\cup T\cup A \cup B$. By the positions of chosen nodes in (1)-(3), there are two pairs $\{b_1^0,t_1\}$ and $\{b_2^0,t_2\}$ in $BH_{n-1}^0$, $\{s_1,a_1^1\}$ and $\{s_3,t_3\}$ in $BH_{n-1}^1$, $\{s_2,a_2^3\}$ and $\{b_1^3,a_1^3\}$ in $BH_{n-1}^3$, respectively. From Lemma~\ref{2-DPC} for $n\geq 3$, $BH_{n-1}^0$ has two disjoint paths, say $W_1[b_1^0,t_1]$ and $W_2[b_2^0,t_2]$, such that $V(W_1[b_1^0,t_1])\cup V(W_2[b_2^0,t_2])=V(BH_{n-1}^0)$; $BH_{n-1}^1$ has two disjoint paths, say $P_1[s_1,a_1^1]$ and $P_3[s_3,t_3]$, such that $V(P_1[s_1,a_1^1])\cup V(P_3[s_3,t_3])=V(BH_{n-1}^1)$; $BH_{n-1}^3$ has two disjoint paths, say $R_2[s_2,a_2^3]$ and $R_1[b_1^3,a_1^3]$, such that $V(R_2[s_2,a_2^3])\cup V(R_1[b_1^3,a_1^3])=V(BH_{n-1}^3)$. There is one pair $b_1^2,a_1^2$ in $BH_{n-1}^2$. By Lemma~\ref{H-l}, for $n\geq 3$, there exists one Hamiltonian path, say $Q_1[b_1^2,a_1^2]$, in $BH_{n-1}^2$. Note that $(a_{1}^{1},b_{1}^2)$, $(a_{1}^{2},b_{1}^3)$, $(a_{1}^{3},b_{1}^0)$ and $(a_{2}^{3},b_{2}^0)$ are crossing edges. Then $BH_{n}$ has a paired three-disjoint path cover relative to $S\cup T$ for $n\geq 3$, in which the three node disjoint paths are as follows:
\\$P[s_1,t_1]=\langle s_1,P_1[s_1,a_1^1],a_1^1,b_{1}^2,Q_1[b_1^2,a_1^2],a_1^2,b_1^3,R_1[b_1^3,a_1^3],a_1^3,b_1^0,W_1[b_1^0,t_1],t_1 \rangle$,
\\$P[s_2,t_2]=\langle s_2,R_2[s_2,a_2^3],a_2^3,b_2^0, W_2[b_2^0,t_2],t_2 \rangle$,
\\$P[s_3,t_3]=\langle s_3,P_3[s_3,t_3],t_3 \rangle$.

{\bf Subcase 1.2.} $f_0=0$, $0<f_3<2$.

It implies that $f_3=1$ and $f_1+f_2=3$. Then $(f_1,f_2)\in\{(0,3),(3,0),(1,2),(2,1)\}$. By the definition of $f_i$, for each $i\in \{0,1,2,3\}$ we have $1\leq \beta_{i}\leq 3$ and there is at least one $\beta_{i}$ equal to 3. The outlet of the proof for each situation is that for $n\geq 4$, by Lemma~\ref{H-l}, Lemma~\ref{2-DPC} and the inductive hypothesis, $BH_{n-1}^i$ has a path cover with $\beta_{i}=1$, $\beta_{i}=2$ and $\beta_{i}=3$ respectively; for $n=3$, by Lemma~\ref{H-l}, Lemma~\ref{2-DPC} and Lemma~\ref{BH2}, $BH_{n-1}^i$ has a path cover with $\beta_{i}=1$, $\beta_{i}=2$ and $\beta_{i}=3$ respectively.

Since the discussions are similar, we only consider one situation that $\beta_{3}=\beta_{0}=2,\beta_{1}=2$ and $\beta_{2}=1$, without loss of generality, let $s_{1}\in V(BH_{n-1}^1),s_{2},s_{3}\in V(BH_{n-1}^3)$ and $t_{1},t_{2}\in V(BH_{n-1}^0),t_{3}\in V(BH_{n-1}^1)$ (see Figure~\ref{Fig-4}).

There are two pairs $\{s_1,a_1^1\}$ and $\{b_3^1,s_3\}$ in $BH_{n-1}^1$. By Lemma~\ref{2-DPC}, for $n\geq 3$, $BH_{n-1}^1$ has two disjoint paths, say $P_1[s_1,a_1^1]$ and $P_3[b_3^1,s_3]$, which contain all the nodes in $BH_{n-1}^1$. There is one pair $\{b_1^2,a_1^2\}$ in $BH_{n-1}^2$, by Lemma~\ref{H-l}, for $n\geq 3$, there exist one Hamiltonian path, say $Q_1[b_1^2,a_1^2]$, in $BH_{n-1}^2$.
There are three pairs $\{b_1^0,t_1\},\{b_2^0,t_2\},\{b_3^0,a_3^0\}$ in $BH_{n-1}^0$, and $\{b_1^3,a_1^3\},\{s_2,a_2^3\},\{s_3,a_3^3\}$ in $BH_{n-1}^3$, respectively. By the inductive hypothesis for $n\geq4$ and by Lemma~\ref{BH2} for $n=3$, $BH_{n-1}^0$ has three disjoint paths, say $W_1[b_1^0,t_1]$, $W_2[b_2^0,t_2]$ and $W_3[b_3^0,a_3^0]$, which contain all the nodes in $BH_{n-1}^0$. Similarly, $BH_{n-1}^3$ has three disjoint paths $R_1[b_1^3,a_1^3]$, $R_2[s_2,a_2^3]$ and $R_3[s_3,a_3^3]$, which contain all the nodes in $BH_{n-1}^0$. Furthermore,  $(a_{3}^{0},b_{3}^1)$, $(a_{1}^{1},b_{1}^2)$, $(a_{1}^{2},b_{1}^3)$, $(a_{1}^{3},b_{1}^0)$, $(a_{2}^{3},b_{2}^0)$ and $(a_{3}^{3},b_{3}^0)$ are crossing edges. Then for $n\geq 3$, $BH_{n}$ has a paired three-disjoint path cover relative to $S\cup T$, in which three node disjoint paths are as follows:
\\$P[s_1,t_1]=\langle s_1,P_1[s_1,a_1^1],a_1^1,b_{1}^2,Q_1[b_1^2,a_1^2],a_1^2,b_1^3,R_1[b_1^3,a_1^3],a_1^3,b_1^0,W_1[b_1^0,t_1],t_1 \rangle$;
\\$P[s_2,t_2]=\langle s_2,R_2[s_2,a_2^3],a_2^3,b_2^0,W_2[b_2^0,t_2],t_2 \rangle$;
\\$P[s_3,t_3]=\langle s_3,R_3[s_3,a_3^3],a_3^3,b_3^0,W_3[b_3^0,a_3^0],a_3^0,b_3^1,P_1[b_3^1,t_3],t_3\rangle$.

{\bf Subcase 1.3.} $f_0>0$, $0\leq f_3<2$.

\begin{figure}
\begin{minipage}[t]{0.5\linewidth}
\centering
\includegraphics[height=5cm,width=7cm]{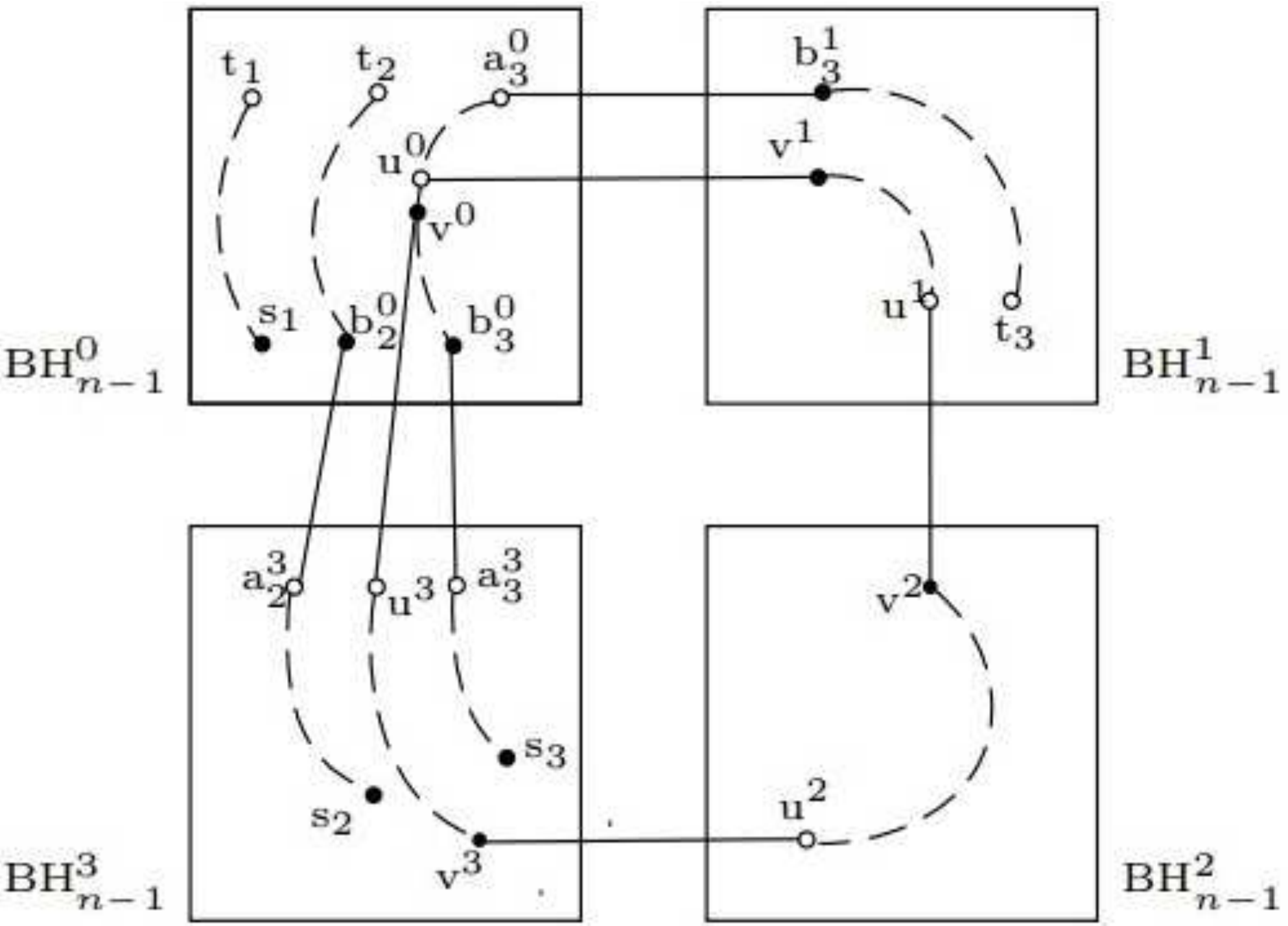}
\vskip-0.3cm
\caption{Illustration of 1.3. in Theorem~\ref{3-DPC}}
\label{Fig-5}
\end{minipage}
\begin{minipage}[t]{0.5\linewidth}
\centering
\includegraphics[height=5cm,width=6.5cm]{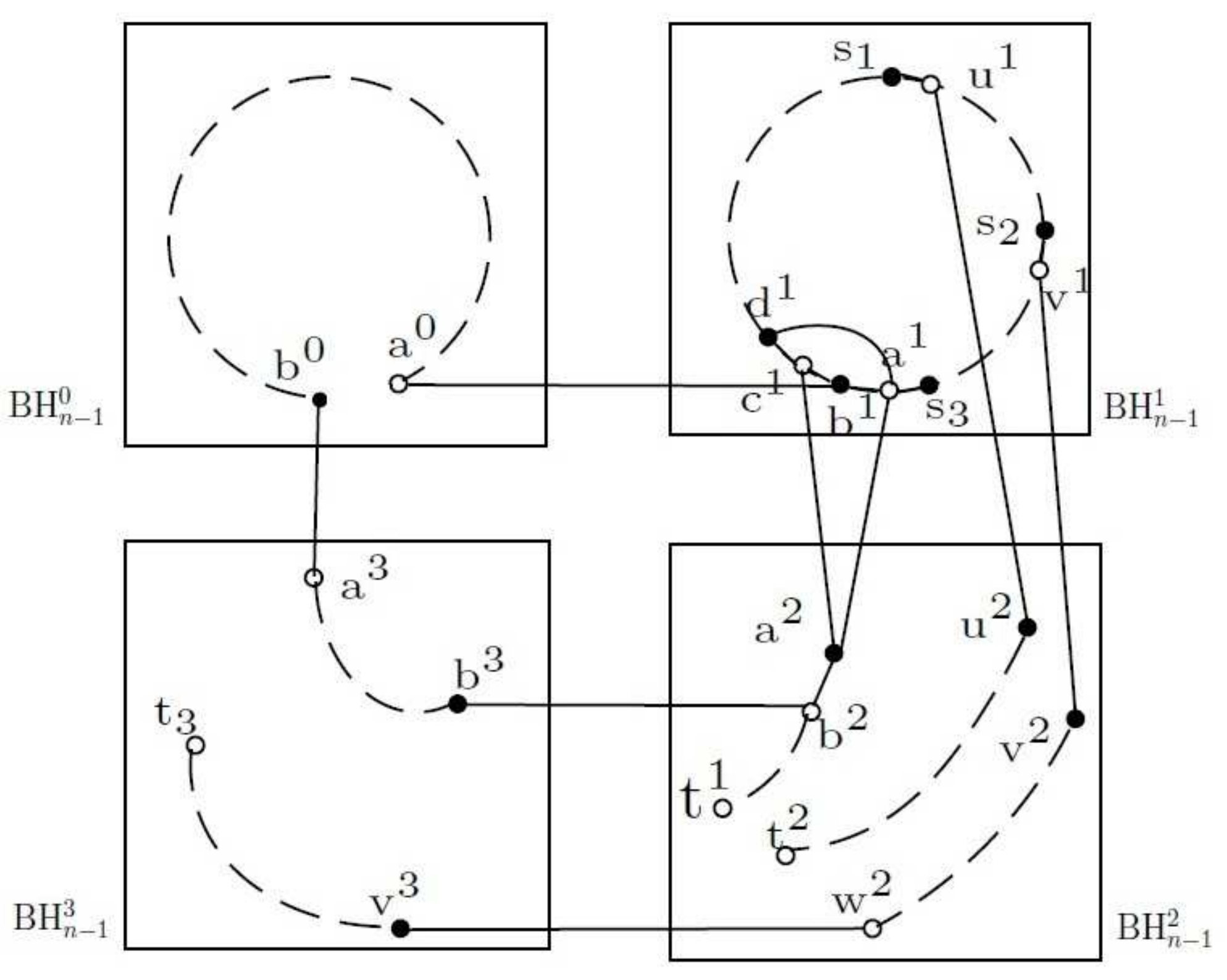}
\vskip-0.3cm
\caption{Illustration of 2.3. in Theorem~\ref{3-DPC}}
\label{Fig-2}
\end{minipage}
\end{figure}

Since $0\leq f_3<2$, we have $f_3=0$ or $1$. There are many situations with different $(f_0,f_1,f_2,f_3)$ in this case. By the definition of $f_i$, there are at least one of $\beta_{i}$ is equal to 0 and no more than one of $\beta_{i}$ is equal to 3, where $i\in \{0,1,2,3\}$. Since the discussions are similar, we only consider one situation that $(f_0,f_1,f_2,f_3)=(1,1,1,1)$, and assume $s_{1}\in V(BH_{n-1}^0),s_{2},s_{3}\in V(BH_{n-1}^3)$ and $t_{1},t_{2}\in V(BH_{n-1}^0),t_{3}\in V(BH_{n-1}^1)$, (see Figure~\ref{Fig-5}).

\f There are three pairs $\{s_1,t_1\}, \{b_2^0,t_2\}$ and $\{b_3^0,a_3^0\}$ in $BH_{n-1}^0$. Similar to the proof of the Subcase 2.2, we can obtain that for $n\geq 3$, $BH_{n-1}^0$ has three disjoint paths, say $W_1[s_1,t_1]$, $W_2[b_2^0,t_2]$ and $W_3[b_3^0,a_3^0]$, which contain all the nodes in $BH_{n-1}^0$. Choosing an edge $(u^0,v^0)$ in $W_3[b_3^0,a_3^0]$, where $u^0$ is white and $v^0$ is black, such that $W_3[b_3^0,a_3^0]=\langle b_3^0,W_3^1[b_3^0,v^0], v^0,u^0,W_3^2[u_0,a_3^0]\rangle$, so there exists one pair of nodes $v^1,u^1$ (resp., $v^3,u^3$) in $BH_{n-1}^1$ (resp., $BH_{n-1}^3$) such that the edge $(v^1,u^0)$ (resp., $(v^0,u^3)$) is a crossing edge. Let $v^2$ (resp., $u^2$) is the neighbor of $u^1$ (resp., $v^3$) in $BH_{n-1}^2$. By Lemma~\ref{H-l}, there exists a Hamiltonian path, say $Q_3[v^2,u^2]$, in $BH_{n-1}^2$ between $v^2$ and $u^2$. Since there are three pairs $\{s_2,a_2^3\},\{v_3,u_3\}$ and $\{s_3,a_3^3\}$ in $BH_{n-1}^3$, similar to the proof of the subcase 1.2, we can obtain that $BH_{n-1}^3$ has three disjoint paths, say $R_2[s_2,a_2^3]$, $R_3^\prime [v_3,u_3]$ and $R_3[s_3,a_3^3]$, which contain all the nodes in $BH_{n-1}^3$ for $n\geq 3$. There are two pairs $\{b_3^1,t_3\}$ and $\{v^1,u^1\}$ in $BH_{n-1}^1$, by Lemma~\ref{2-DPC}, $BH_{n-1}^1$ has two disjoint paths, say $P_3[b_3^1,t_3]$ and $P_3^\prime[v^1,u^1]$, which contain all the nodes in $BH_{n-1}^1$. Furthermore, $(a_{3}^{0},b_{3}^1),(u^{0},v^1),(u^{1},v^2),(u^{2},v^3),(u^{3},v^0),(a_{3}^{3},b_{3}^0)$ and $(a_{2}^{3},b_{2}^0)$ are crossing edges. Then  $BH_{n}$ for $n\geq 3$ has a paired three-disjoint path cover relative to $S\cup T$ in which three node disjoint paths are as follows:
\\$P[s_1,t_1]=W_1[s_1,t_1]$,
\\$P[s_2,t_2]=\langle s_2,R_2[s_2,a_2^3],a_2^3,b_2^0,W_2[b_2^0,t_2],t_2 \rangle$;
\\$P[s_3,t_3]=\langle s_3,R_3[s_3,a_3^3],a_3^3,b_{3}^0,W_3^1[b_3^0,v^0],v^0,u^3,R_3^\prime[u^3,v^3],v^3,u^2,Q_3[u^2,v^2],v^2,u^1,P_3^\prime[u^1,v^1],
\\v^1,u^0,W_3^2[u_0,a_3^0],a_3^0,b_3^1,P_3[b_3^1,t_3],t_3\rangle$.

\medskip
{\bf Case 2.} $f_0\geq 1$ and $f_3\geq 2$.

Since $0\leq f_k\leq 4$ for $k\in \{0,1,2,3\}$ and $\sum\limits_{i=0}^{3} f_i=4$, we have $f_0=f_3=2$ or $f_0=1,f_3=3$ or $f_0=1,f_3=2$. Now let $s_j\in V(BH_{n-1}^{g_j}),t_j\in V(BH_{n-1}^{h_j})$ for $j\in\{1,2,3\}$, where $g_j,h_j\in\{0,1,2,3\}$.

{\bf Subcase 2.1.} $f_0=f_3=2$.

It implies that there are two $BH_{n-1}^i$'s with $\beta_{i}=0$ and two $BH_{n-1}^i$'s with $\beta_{i}=3$, $i\in \{0,1,2,3\}$. If the two $BH_{n-1}^i$'s with $\beta_{i}=0$ are adjacent, without loss of generality, we can assume $\beta_0=\beta_1=0$, then $\beta_2=\beta_3=3$. Since $\beta_{i}=M_{1}^i+M_{2}^i+M_{3}^i$, then $M_{j}^0=M_{j}^1=0$ and $M_{j}^2=M_{j}^3=1$, $j\in \{1,2,3\}$. By the definition of $M_{j}^i$, for each $j\in \{1,2,3\}$, $\{0,1\}\subseteq \{h_j+1,h_j+2,\ldots,g_j-1\}$ and $\{2,3\}\subseteq \{g_j,g_j+1,\ldots,h_j\}$. Since there are four elements in the set $\{g_j,g_j+1,\ldots,h_j,h_j+1,\ldots,g_j-1\}$, then we have $g_1=g_2=g_3=2$ and $h_1=h_2=h_3=3$. Thus the three white nodes $t_{1},t_{2},t_{3}$ are all in $BH_{n-1}^3$, which contradicts with the assumption that $t_{1},t_{2},t_{3}$ are not in the same $BH_{n-1}^i$ for $i\in\{0,1,2,3\}$.

If the two $BH_{n-1}^i$'s with $\beta_{i}=0$ are not adjacent, without loss of generality, we can assume $\beta_0=\beta_2=0$, then $\beta_1=\beta_3=3$. Since $\beta_{i}=M_{1}^i+M_{2}^i+M_{3}^i$, then $M_{j}^0=M_{j}^2=0$ and $M_{j}^1=M_{j}^3=1$ for each $j\in \{1,2,3\}$. By the definition of $M_{j}^i$,  for each $j\in \{1,2,3\}$, $\{0,2\}\subseteq \{h_j+1,h_j+2,\ldots,g_j-1\}$ and $\{1,3\}\subseteq \{g_j,g_j+1,\ldots,h_j\}$. Since the sequence $g_j,g_j+1,\ldots,h_j$ are continuous, then we can deduce $2$ or $0\in\{g_j,g_j+1,\ldots,h_j\}$ because of $\{1,3\}\subseteq \{g_j,g_j+1,\ldots,h_j\}$, which contradicts with $\{0,2\}\subseteq \{h_j+1,h_j+2,\ldots,g_j-1\}$.

{\bf Subcase 2.2.} $f_0=1,f_3=3$.

It implies that there are one $BH_{n-1}^i$ with $\beta_{i}=0$ and three $BH_{n-1}^i$'s with $\beta_{i}=3$, $i\in \{0,1,2,3\}$. Without loss of generality, we can assume $\beta_0=0$, then $\beta_1=\beta_2=\beta_3=3$. Since $\beta_{i}=M_{1}^i+M_{2}^i+M_{3}^i$, then $M_{j}^0=0$ and $M_{j}^1=M_{j}^2=M_{j}^3=1$, $j\in \{1,2,3\}$. By the definition of $M_{j}^i$, for each $j\in \{1,2,3\}$, we have $0\in \{h_j+1,h_j+2,\ldots,g_j-1\}$ and $\{1,2,3\}\subseteq \{g_j,g_j+1,\ldots,h_j\}$. Since there are four elements in the set $\{g_j,g_j+1,\ldots,h_j,h_j+1,\ldots,g_j-1\}$, then $g_1=g_2=g_3=1$ and $h_1=h_2=h_3=3$. Thus the three white nodes $t_{1},t_{2},t_{3}$ are all in $BH_{n-1}^3$, which contradicts with the assumption that $t_{1},t_{2},t_{3}$ are not in the same $BH_{n-1}^i$ for $i\in\{0,1,2,3\}$.

{\bf Subcase 2.3.} $f_0=1,f_3=2$.

It implies that there are one $BH_{n-1}^i$ with $\beta_{i}=0$ and two $BH_{n-1}^i$'s with $\beta_{i}=3$, $i\in \{0,1,2,3\}$. Without loss of generality, we may assume $\beta_{0}=0 $, i.e., $M_{j}^0=0$, $j\in \{1,2,3\}$. By the definition of $M_{j}^0$, for each $j\in \{1,2,3\}$, we have $0\in \{h_j+1,h_j+2,\ldots,g_j-1\}$. If $\beta_{3}=3 $, similarly, for each $j\in \{1,2,3\}$, we have $3\in \{g_j,g_j+1,\ldots,h_j\}$. If $h_j=1$, $h_j+1=2$, then $\{0,2\}\in \{h_j+1,h_j+2,\ldots,g_j-1\}$ for each $j\in \{1,2,3\}$. Since the the sequence $h_j+1,h_j+2,\ldots,g_j-1$ are continuous, and $3\in \{g_j,g_j+1,\ldots,h_j\}$, then we can deduce that for each $j\in \{1,2,3\}$, $\{0,1,2\}\in\{h_j+1,h_j+2,\ldots,g_j-1\}$. Thus $M_{j}^0=M_{j}^1=M_{j}^2=0$ and $\beta_{0}=\beta_{1}=\beta_{2}=0$, i.e., $f_0=3$ which is a contradiction. If $h_j=2$, by the similar discussion as that of $h_j=1$, we can obtain a contradiction. Thus $h_j=3$ for each $j\in \{1,2,3\}$, i.e., the three white nodes $t_{1},t_{2},t_{3}$ are all in $BH_{n-1}^3$, which contradicts with the assumption that $t_{1},t_{2},t_{3}$ are not in the same $BH_{n-1}^i$ for $i\in\{0,1,2,3\}$. Thus, we have $\beta_{3}\neq3 $, then $\beta_{1}=\beta_{2}=3$.

Now we only need to consider $\beta_{0}=0$ and $\beta_{1}=\beta_{2}=3$. By the definitions of $\beta_{i}$ and $M_{j}^i$, for each $j\in \{1,2,3\}$, $g_j=1,h_j\geq 2$, i.e., $s_1,s_2$ and $s_3$ are all in $BH_{n-1}^1$ and $t_1,t_2$ and $t_3$ are in $BH_{n-1}^2$ or $BH_{n-1}^3$. Since the discussion of situations are similar, assume $s_{1},s_{2},s_{3}\in V(BH_{n-1}^1)$ and $t_{1},t_{2}\in V(BH_{n-1}^2)$, $t_{3}\in V(BH_{n-1}^3)$, see Figure~\ref{Fig-2}.

By Lemma~\ref{5-path}, there exists a Hamiltonian cycle $C_{1}$ contains the path, say $\langle s_j,a^1,b^1,c^1,d^1\rangle$ in $BH_{n-1}^{1}$, where $\langle a^1,b^1,c^1,d^1\rangle$ is a $4$-cycle, $a^1$ (resp. $b^1$) is symmetric with $c^1$ (resp. $d^1$) and $b^1,d^1$ are not both in the set $\{s_1,s_2,s_3\}$. Without loss of generality, we may assume the location of $s_1,s_2$ is as the case in Figure~\ref{Fig-2} and $b^1\notin \{s_1,s_2\}$ ($d^1$ maybe equal to $s_{1}$).
Let $a^0$ in $BH_{n-1}^0$ be a neighbor of $b^1$. Let $a^2\in V(BH_{n-1}^2)$ be a common neighbor of $a^1$ and $c^1$. Note that $a^1$ and $d^1$ are adjacent. Let $u^1$,$v^1$ be the neighbors of $s_{1}$ and $s_{2}$ in  $C_{1}$, respectively. Let $C_{1}=\langle a^1,b^1,c^1,d^1,P_1[d^1,s_{1}],s_{1},u^1,P_2[u^1,s_{2}],s_{2},v^1,P_3[v^1,s_{3}],s_{3},a^1\rangle$.
Let $u^2,v^2\in V(BH_{n-1}^2)$ be the neighbors of $u^1$ and $v^1$, respectively. Let $w^2\in BH_{n-1}^2$ be a white node, such that $w^2\notin\{t_{1},t_2\}$. For $n\geq 4$, by the inductive hypothesis (for $n=3$, by Lemma~\ref{BH2} respectively), $BH_{n-1}^2$ has three disjoint paths, say $Q_1[a^2,t_{1}]$, $Q_2[u^2,t_{2}]$ and $Q_3[v^2,w^2]$ which contain all the nodes in $BH_{n-1}^2$. Assume that $b^2$ ($b^2$ may be equal to $t_{1}$) in $Q_1[a^2,t_{1}]$ is a neighbor of $a^2$. Then $Q_1[a^2,t_{1}]=\langle a^2,b^2,Q_1[b^2,t_{1}],t_{1}\rangle$. Let $b^3,v^3\in BH_{n-1}^3$ such that the edges $(b^2,b^3),(w^2,v^3)$ are crossing edges. Assume that $a^3$ is a white node in $BH_{n-1}^3$ and $a^3\neq t_{3}$, then by Lemma~\ref{2-DPC}, $BH_{n-1}^3$ has two disjoint paths, say $R_1[b^3,a^3]$, $R_3[v^3,t_3]$ which contain all the nodes in $BH_{n-1}^3$. Let $b^0\in BH_{n-1}^0$ be a neighbor of $a^3$ and $a^0\in BH_{n-1}^0$ be a neighbor of $b^1$, by Lemma~\ref{H-l}, there is a Hamiltonian path $W_1[a^0,b^0]$ between $a^0$ and $b^0$.

Then for $n\geq3$, $BH_{n}$ has a paired three-disjoint path cover relative to $S\cup T$ in which three node disjoint paths are as follows:
\\$P[s_1,t_1]=\langle P_1[s_1,d^1],d^1,a^1,a^2,c^1,b^1,a^0,W_1[a^0,b^0],b^0,a^3,R_1[a^3,b^3],b^3,c^2,Q_1[c^2,t_1],t_1 \rangle$,
\\$P[s_2,t_2]=\langle P_1[s_2,u^1],u^1,Q_2[u^2,t_2],t_2 \rangle$,
\\$P[s_3,t_3]=\langle P_3[s_3,v^1],v^1,v^2,Q_3[v^2,w^2],w^2,v^3,R_3[v^3,t_3],t_3 \rangle$.

By the above cases, the proof is complete. \hfill\qed

\section{Conclusion}

In this paper, the paired three-disjoint path cover of $BH_{n}$ is obtained. The result shows that for $n\geq 3$, let $S=\{s_{1}, s_{2}, s_{3}\}\subseteq B$ and $T=\{t_{1}, t_{2}, t_{3}\}\subseteq W$, then $BH_{n}$ has a paired 3-disjoint path cover relative to $S\cup T$, where $B$ and $W$ are two distinct partite sets of $BH_{n}$. If $t_2$ is adjacent to $s_3$, then we get the paired 2-disjoint path cover; if $t_1$ adjacent to $s_2$ and $t_2$ is adjacent to $s_3$, then we get the paired 1-disjoint path cover, so the result in this paper is more general than the known results about the paired $k$-disjoint path covers of $BH_{n}$ for $n\geq 3$ and $k=1,2$.

On the other hand, because the paired 3-disjoint path covers of $BH_{n}$ for $n=2$ relative to some sets $S\cup T$ are not exist, the low bound of $n\geq 3$ in the result is the best. In fact, as an example, in $BH_2$, let $S=\{s_1, s_2, s_3\}$ and $T=\{t_1, t_2, t_3\}$, where $s_1=(1,0), s_2=(3,0), s_3=(1,2)$, $t_1=(0,0), t_2=(2,0)$ and $t_3=(0,1)$. For $i\in\{1,2,3\}$, let $P[t_i,s_i]=P_i$. If there exists a paired 3-disjoint path cover $\{P_1,P_2,P_3\}$ relative to $S\cup T$, then at least one of them, say $P_i$, satisfies $|P_i|>1$. If $|P_1|>1$, by these three paths being node-disjoint and symmetricity of $BH_2$, $P_1=\langle t_1, (1,1), (2,1), (3,2), \ldots, s_1\rangle$. It implies that $P_2=\langle t_2,s_2\rangle$ and $P_3=\langle t_3,s_3\rangle$. But there does not exist a path $P_1$ which contains all the nodes in $BH_2$ except $\{t_2,s_2,t_3,s_3\}$. The discussion for the case $|P_2|>1$ or $|P_3|>1$ is similar. As a result, there does not exist a paired 3-disjoint path cover in $BH_{2}$ relative to $S\cup T$.

Further more, the paired $k$-disjoint path cover of $BH_{n}$ for $k\geq 4, n\geq k$ need to be studied in the future.



\section*{Appendix A. All the cases in the proof of Lemma 2.9}

\texttt{For convenience, let $a=(0,0)$, $b=(1,0)$, $c=(2,0)$, $d=(3,0)$, $e=(0,1)$, $f=(1,1)$, $g=(2,1)$, $h=(3,1)$, $i=(0,2)$, $j=(1,2)$, $k=(2,2)$, $l=(3,2)$, $m=(0,3)$, $n=(1,3)$, $o=(2,3)$ and $p=(3,3)$.}

{\small
\begin{center}
\texttt{Table 1. $t_1=(0,0)$ and $t_2=(2,0)$}
\medskip

\begin{tabular}{|l|l|l|l|l|l|l|}

\hline
$s_1$ & $s_2$ & $s_3$ & $t_3$ &$P_1[s_1,t_1]$ & $P_2[s_2,t_2]$ & $P_3[s_3,t_3]$ \\
\hline
$b$ & $d$ & $f$ & $e$ & $\langle b,a\rangle$ & $ \langle d,m,p,o,n,i,l,k,j,g,h,c\rangle$ & $\langle f,e\rangle$ \\
\hline
$b$ & $d$ & $j$ & $k$ & $\langle b,a\rangle$ & $\langle d,m,p,o,n,i,l,g,f,e,h,c\rangle$ & $\langle j,k\rangle$ \\
\hline
$b$ & $d$ & $n$ & $o$ & $\langle b,a\rangle$ & $\langle d,m,p,k,l,i,j,g,f,e,h,c\rangle$ & $\langle n,o\rangle$ \\
\hline
$f$ & $h$ & $b$ & $e$ & $\langle f,g,l,i,p,m,d,a\rangle$ & $\langle h,c\rangle$ & $\langle b,o,n,k,j,e\rangle$ \\
\hline
$f$ & $h$ & $j$ & $g$ & $\langle f,e,l,i,n,k,p,m,b,o,d,a\rangle$ & $\langle h,c\rangle$ & $\langle j,g\rangle$ \\
\hline
$f$ & $h$ & $n$ & $k$ & $\langle f,e,j,g,l,i,p,o,d,m,b,a\rangle$ & $\langle h,c\rangle$ & $\langle n,k\rangle$ \\
\hline
$j$ & $l$ & $b$ & $o$ & $\langle j,g,f,e,h,a\rangle$ & $ \langle l,k,n,i,p,m,d,c\rangle$ & $\langle b,o\rangle$ \\
\hline
$j$ & $l$ & $f$ & $e$ & $\langle j,k,p,m,d,a\rangle$ & $ \langle l,i,n,o,b,c\rangle$ & $\langle f,g,h,e\rangle$ \\
\hline
$j$ & $l$ & $n$ & $k$ & $\langle j,g,h,e,f,a\rangle$ & $ \langle l,i,p,o,d,m,b,c\rangle$ & $\langle n,k\rangle$ \\
\hline
$n$ & $p$ & $b$ & $m$ & $\langle n,o,d,a\rangle$ & $ \langle p,i,l,k,j,g,f,e,h,c\rangle$ & $\langle b,m\rangle$ \\
\hline
$n$ & $p$ & $f$ & $k$ & $\langle n,m,d,o,b,a\rangle$ & $ \langle p,i,j,g,h,c\rangle$ & $\langle f,e,l,k\rangle$ \\
\hline
$n$ & $p$ & $j$ & $k$ & $\langle n,o,b,m,d,a\rangle$ & $ \langle p,i,l,g,h,e,f,c\rangle$ & $\langle j,k\rangle$ \\
\hline
$b$ & $h$ & $d$ & $o$ & $\langle b,m,p,i,n,k,l,g,j,e,f,a\rangle$ & $ \langle h,c\rangle$ & $\langle d,o\rangle$ \\
\hline
$b$ & $h$ & $f$ & $e$ & $\langle b,a\rangle$ & $ \langle h,g,j,k,l,i,n,o,p,m,d,c\rangle$ & $\langle f,e\rangle$ \\
\hline
$b$ & $h$ & $j$ & $k$ & $\langle b,a\rangle$ & $ \langle h,e,f,g,l,i,n,o,p,m,d,c\rangle$ & $\langle j,k\rangle$ \\
\hline
$b$ & $h$ & $n$ & $o$ & $\langle b,a\rangle$ & $ \langle h,e,f,g,j,k,l,i,p,m,d,c\rangle$ & $\langle n,o\rangle$ \\
\hline
$f$ & $d$ & $b$ & $o$ & $\langle f,a\rangle$ & $ \langle d,m,p,i,n,k,l,g,j,e,h,c\rangle$ & $\langle b,o\rangle$ \\
\hline
$f$ & $d$ & $h$ & $i$ & $\langle f,a\rangle$ & $ \langle d,m,p,o,b,c\rangle$ & $\langle h,e,l,g,j,k,n,i\rangle$ \\
\hline
$f$ & $d$ & $j$ & $e$ & $\langle f,g,h,a\rangle$ & $ \langle d,m,p,o,b,c\rangle$ & $\langle j,i,n,k,l,e\rangle$ \\
\hline
$f$ & $d$ & $n$ & $i$ & $\langle f,a\rangle$ & $ \langle d,m,b,o,p,k,l,g,j,e,h,c\rangle$ & $\langle n,i\rangle$ \\
\hline
$b$ & $l$ & $d$ & $m$ & $\langle b,o,p,i,n,k,j,e,f,a\rangle$ & $ \langle l,g,h,c\rangle$ & $\langle d,m\rangle$ \\
\hline
$b$ & $l$ & $f$ & $g$ & $\langle b,a\rangle$ & $ \langle l,k,j,i,n,m,p,o,d,c\rangle$ & $\langle f,e,h,g\rangle$ \\
\hline
$b$ & $l$ & $j$ & $k$ & $\langle b,o,n,i,p,m,d,a\rangle$ & $ \langle l,g,f,e,h,c\rangle$ & $\langle j,k\rangle$ \\
\hline
$b$ & $l$ & $n$ & $k$ & $\langle b,o,p,m,d,a\rangle$ & $ \langle l,i,j,g,h,e,f,c\rangle$ & $\langle n,k\rangle$ \\
\hline
$j$ & $d$ & $b$ & $i$ & $\langle j,k,l,e,f,g,h,a\rangle$ & $ \langle d,c\rangle$ & $\langle b,m,p,o,n,i\rangle$ \\
\hline
$j$ & $d$ & $f$ & $g$ & $\langle j,k,l,i,n,o,p,m,b,a\rangle$ & $ \langle d,c\rangle$ & $\langle f,e,h,g\rangle$ \\
\hline
$j$ & $d$ & $l$ & $o$ & $\langle j,e,h,g,f,a\rangle$ & $ \langle d,m,b,c\rangle$ & $\langle l,k,n,i,p,o\rangle$ \\
\hline
$j$ & $d$ & $n$ & $o$ & $\langle j,k,l,e,f,g,h,a\rangle$ & $ \langle d,c\rangle$ & $\langle n,i,p,m,b,o\rangle$ \\
\hline
$b$ & $p$ & $d$ & $m$ & $\langle b,a\rangle$ & $ \langle p,o,n,i,j,k,l,g,h,e,f,c\rangle$ & $\langle d,m\rangle$ \\
\hline
\end{tabular}
\end{center}

\begin{center}
\begin{tabular}{|l|l|l|l|l|l|l|}
\hline
$s_1$ & $s_2$ & $s_3$ & $t_3$ & $P_1[s_1,t_1]$ & $P_2[s_2,t_2]$ & $P_3[s_3,t_3]$ \\
\hline
$b$ & $p$ & $f$ & $g$ & $\langle b,a\rangle$ & $ \langle p,m,d,o,n,i,l,k,j,e,h,c\rangle$ & $\langle f,g\rangle$ \\
\hline
$b$ & $p$ & $j$ & $i$ & $\langle b,a\rangle$ & $ \langle p,m,d,o,n,k,l,e,h,g,f,c\rangle$ & $\langle j,i\rangle$ \\
\hline
$b$ & $p$ & $n$ & $m$ & $\langle b,a\rangle$ & $ \langle p,i,l,k,j,g,h,e,f,c\rangle$ & $\langle n,o,d,m\rangle$ \\
\hline
$n$ & $d$ & $b$ & $m$ & $\langle n,o,p,i,j,k,l,e,h,g,f,a\rangle$ & $ \langle d,c\rangle$ & $\langle b,m\rangle$ \\
\hline
$n$ & $d$ & $f$ & $e$ & $\langle n,o,p,i,l,k,j,g,h,a\rangle$ & $ \langle d,m,b,c\rangle$ & $\langle f,e\rangle$ \\
\hline
$n$ & $d$ & $j$ & $k$ & $\langle n,o,p,i,l,e,f,g,h,a\rangle$ & $ \langle d,m,b,c\rangle$ & $\langle j,k\rangle$ \\
\hline
$n$ & $d$ & $p$ & $o$ & $\langle n,i,j,k,l,e,f,g,h,a\rangle$ & $ \langle d,m,b,c\rangle$ & $\langle p,o\rangle$ \\
\hline
$f$ & $l$ & $b$ & $o$ & $\langle f,g,j,e,h,a\rangle$ & $ \langle l,k,n,i,p,m,d,c\rangle$ & $\langle b,o\rangle$ \\
\hline
$f$ & $l$ & $h$ & $e$ & $\langle f,g,j,i,p,m,d,a\rangle$ & $ \langle l,k,n,o,b,c\rangle$ & $\langle h,e\rangle$ \\
\hline
$f$ & $l$ & $j$ & $e$ & $\langle f,a\rangle$ & $ \langle l,k,n,i,p,o,d,m,b,c\rangle$ & $\langle j,g,h,e\rangle$ \\
\hline
$f$ & $l$ & $n$ & $k$ & $\langle f,e,h,a\rangle$ & $ \langle l,f,j,i,p,o,d,m,b,c\rangle$ & $\langle n,k\rangle$ \\
\hline
$j$ & $h$ & $b$ & $o$ & $\langle j,g,f,a\rangle$ & $ \langle h,e,l,k,n,i,p,m,d,c\rangle$ & $\langle b,o\rangle$ \\
\hline
$j$ & $h$ & $f$ & $g$ & $\langle j,e,l,k,n,i,p,o,b,m,d,a\rangle$ & $ \langle h,c\rangle$ & $\langle f,g\rangle$ \\
\hline
$j$ & $h$ & $l$ & $e$ & $\langle j,k,p,i,n,o,d,m,b,a\rangle$ & $ \langle h,g,f,c\rangle$ & $\langle l,e\rangle$ \\
\hline
$j$ & $h$ & $n$ & $e$ & $\langle j,i,p,o,d,m,b,a\rangle$ & $ \langle h,g,f,c\rangle$ & $\langle n,k,l,e\rangle$ \\
\hline
$f$ & $p$ & $b$ & $m$ & $\langle f,g,h,a\rangle$ & $ \langle p,i,l,e,j,k,n,o,d,c\rangle$ & $\langle b,m\rangle$ \\
\hline
$f$ & $p$ & $h$ & $e$ & $\langle f,a\rangle$ & $ \langle p,i,l,g,j,k,n,o,d,m,b,c\rangle$ & $\langle h,e\rangle$ \\
\hline
$f$ & $p$ & $j$ & $e$ & $\langle f,a\rangle$ & $ \langle p,i,l,k,n,o,d,m,b,c\rangle$ & $\langle j,g,h,e\rangle$ \\
\hline
$f$ & $p$ & $n$ & $k$ & $\langle f,g,l,e,h,a\rangle$ & $ \langle p,o,d,m,b,c\rangle$ & $\langle n,i,j,k\rangle$ \\
\hline
$n$ & $h$ & $b$ & $i$ & $\langle n,k,l,g,j,e,f,a\rangle$ & $ \langle h,c\rangle$ & $\langle b,m,d,o,p,i\rangle$ \\
\hline
$n$ & $h$ & $f$ & $e$ & $\langle n,m,d,a\rangle$ & $ \langle h,g,j,i,l,k,p,o,b,c\rangle$ & $\langle f,e\rangle$ \\
\hline
$n$ & $h$ & $j$ & $k$ & $\langle n,i,p,o,d,m,b,a\rangle$ & $ \langle h,g,f,c\rangle$ & $\langle j,e,l,k\rangle$ \\
\hline
$n$ & $h$ & $p$ & $i$ & $\langle n,m,d,o,b,a\rangle$ & $ \langle h,g,j,k,l,e,f,c\rangle$ & $\langle p,i\rangle$ \\
\hline
$j$ & $p$ & $b$ & $m$ & $\langle j,g,f,a\rangle$ & $ \langle p,i,n,k,l,e,h,c\rangle$ & $\langle b,o,d,m\rangle$ \\
\hline
$j$ & $p$ & $f$ & $e$ & $\langle j,g,h,a\rangle$ & $ \langle p,i,l,k,n,o,d,m,b,c\rangle$ & $\langle f,e\rangle$ \\
\hline
$j$ & $p$ & $l$ & $k$ & $\langle j,g,h,e,f,a\rangle$ & $ \langle p,i,n,o,d,m,b,c\rangle$ & $\langle l,k\rangle$ \\
\hline
$j$ & $p$ & $n$ & $i$ & $\langle j,k,l,g,h,e,f,a\rangle$ & $ \langle p,o,d,m,b,c\rangle$ & $\langle n,i\rangle$ \\
\hline
$n$ & $l$ & $b$ & $m$ & $\langle n,k,j,e,h,g,f,a\rangle$ & $ \langle l,i,p,o,d,c\rangle$ & $\langle b,m\rangle$ \\
\hline
$n$ & $l$ & $f$ & $e$ & $\langle n,k,j,g,h,a\rangle$ & $ \langle l,i,p,o,d,m,b,c\rangle$ & $\langle f,e\rangle$ \\
\hline
$n$ & $l$ & $j$ & $k$ & $\langle n,i,p,o,b,m,d,a\rangle$ & $ \langle l,e,h,g,f,c\rangle$ & $\langle j,k\rangle$ \\
\hline
$n$ & $l$ & $p$ & $i$ & $\langle n,o,d,m,b,a\rangle$ & $ \langle l,k,j,g,h,e,f,c\rangle$ & $\langle p,i\rangle$ \\
\hline
\end{tabular}
\end{center}

\newpage
\begin{center}
\texttt{Table 2. $t_1=(0,0)$ and $t_2=(0,1)$}

\begin{tabular}{|l|l|l|l|l|l|l|}
\hline
$s_1$ & $s_2$ & $s_3$ & $t_3$ &$P_1[s_1,t_1]$ & $P_2[s_2,t_2]$ & $P_3[s_3,t_3]$ \\
\hline
$b$&$d$& $f$&$c$& $\langle b,a\rangle$ & $\langle d,m,p,o,n,i,l,k,j,e\rangle$ &$\langle f,g,h,c\rangle$ \\
\hline
$b$&$d$& $j$ & $i$ & $\langle b,c,h,g,f,a\rangle$ & $\langle d,m,p,o,n,k,l,e\rangle$ & $\langle j,i\rangle$ \\
\hline
$b$&$d$& $n$ & $o$ & $\langle b,c,f,a\rangle$ & $\langle d,m,p,i,l,k,j,h,h,e\rangle$ & $\langle n,o\rangle$ \\
\hline
$f$& $h$& $b$ & $c$ & $\langle f,g,j,k,l,i,n,o,p,m,d,a\rangle$ & $\langle h,e\rangle$ & $\langle b,c\rangle$ \\
\hline
$f$& $h$& $j$ & $k$ & $\langle f,g,l,i,n,o,p,m,d,c,b,a\rangle$ & $ \langle h,e\rangle$ & $\langle j,k\rangle$ \\
\hline
$f$& $h$& $n$ & $o$ & $\langle f,g,j,k,l,i,p,m,d,c,b,a\rangle$ & $\langle h,e\rangle$ & $\langle n,o\rangle$ \\
\hline
$j$& $l$& $b$ & $m$ & $\langle j,g,f,a\rangle$& $ \langle l,k,n,i,p,o,d,c,h,e\rangle$& $\langle b,m\rangle$ \\
\hline
$j$& $l$& $f$ & $g$ & $\langle j,k,n,i,p,o,d,m,b,c,h,a\rangle$& $ \langle l,e\rangle$& $\langle f,g\rangle$ \\
\hline
$j$& $l$& $n$ & $k$ & $\langle j,g,f,a\rangle$& $ \langle l,i,p,o,d,m,b,c,h,e\rangle$& $\langle n,k\rangle$ \\
\hline
$n$& $p$& $b$ & $m$ & $\langle n,o,d,a\rangle$& $ \langle p,i,l,k,j,g,h,c,f,e\rangle$& $\langle b,m\rangle$ \\
\hline
$n$& $p$& $f$ & $c$ & $\langle n,o,d,m,b,a\rangle$& $ \langle p,i,l,k,j,g,h,e\rangle$& $\langle f,c\rangle$ \\
\hline
$n$& $p$& $j$ & $k$ & $\langle n,o,d,m,b,c,f,g,h,a\rangle$& $ \langle p,i,l,e\rangle$& $\langle j,k\rangle$ \\
\hline
$f$& $d$& $b$ & $m$ & $\langle f,c,h,a\rangle$& $ \langle d,o,p,i,n,k,l,g,j,e\rangle$& $\langle b,m\rangle$ \\
\hline
$f$& $d$& $h$ & $m$ & $\langle f,a\rangle$& $ \langle d,o,p,i,n,k,l,g,j,e\rangle$& $\langle h,c,b,m\rangle$ \\
\hline
$f$& $d$& $j$ & $k$ & $\langle f,a\rangle$& $ \langle d,c,b,m,p,o,n,i,l,g,h,e\rangle$& $\langle j,k\rangle$ \\
\hline
$f$& $d$& $n$ & $o$ & $\langle f,a\rangle$& $ \langle d,c,b,m,p,i,j,k,l,g,h,e\rangle$& $\langle n,o\rangle$ \\
\hline
$b$& $h$& $d$ & $o$ & $\langle b,c,f,a\rangle$& $ \langle h,g,l,i,p,m,n,k,j,e\rangle$& $\langle d,o\rangle$ \\
\hline
$b$& $h$& $f$ & $c$ & $\langle b,o,n,k,l,g,j,i,p,m,d,a\rangle$& $ \langle h,e\rangle$& $\langle f,c\rangle$ \\
\hline
$b$& $h$& $j$ & $k$ & $\langle b,c,d,m,p,o,n,i,l,g,f,a\rangle$& $ \langle h,e\rangle$& $\langle j,k\rangle$ \\
\hline
$b$& $h$& $n$ & $o$ & $\langle b,c,d,m,p,i,l,k,j,g,f,a\rangle$& $ \langle h,e\rangle$& $\langle n,o\rangle$ \\
\hline
$j$& $d$& $b$ & $m$ & $\langle j,g,f,c,h,a\rangle$& $ \langle d,o,p,i,n,k,l,e\rangle$& $\langle b,m\rangle$ \\
\hline
$j$& $d$& $f$ & $o$ & $\langle j,g,h,a\rangle$& $ \langle d,m,p,i,n,k,l,e\rangle$& $\langle f,c,b,o\rangle$ \\
\hline
$j$ & $d$ & $l$ & $k$ & $\langle j,g,h,a\rangle$ & $ \langle d,m,n,i,p,o,b,c,f,e\rangle$ & $\langle l,k\rangle$ \\
\hline
$j$ & $d$ & $n$ & $i$ & $\langle j,g,f,a\rangle$ & $ \langle d,m,p,o,b,c,f,e\rangle$ & $\langle n,k,l,i\rangle$ \\
\hline
$b$ & $l$ & $d$ & $o$ & $\langle b,a\rangle$ & $ \langle l,i,p,m,n,k,j,g,h,c,f,e\rangle$ & $\langle d,o\rangle$ \\
\hline
$b$ & $l$ & $f$ & $c$ & $\langle b,o,n,k,j,i,p,m,d,a\rangle$ & $ \langle l,g,h,e\rangle$ & $\langle f,c\rangle$ \\
\hline
$b$ & $l$ & $j$ & $k$ & $\langle b,o,n,i,p,m,d,c,h,a\rangle$ & $ \langle l,g,f,e\rangle$ & $\langle j,k\rangle$ \\
\hline
$b$ & $l$ & $n$ & $i$ & $\langle b,o,p,m,d,c,h,a\rangle$ & $ \langle l,k,j,g,f,e\rangle$ & $\langle n,i\rangle$ \\
\hline
$n$ & $d$ & $b$ & $m$ & $\langle n,o,p,i,l,k,j,g,f,a\rangle$ & $ \langle d,c,h,e\rangle$ & $\langle b,m\rangle$ \\
\hline
$n$ & $d$ & $f$ & $m$ & $\langle n,i,j,g,h,a\rangle$ & $ \langle d,o,p,k,l,e\rangle$ & $\langle f,c,b,m\rangle$ \\
\hline
$n$ & $d$ & $j$ & $k$ & $\langle n,i,l,g,f,a\rangle$ & $ \langle d,o,p,m,b,c,h,e\rangle$ & $\langle j,k\rangle$ \\
\hline
$n$ & $d$ & $p$ & $o$ & $\langle n,i,l,k,j,g,f,a\rangle$ & $ \langle d,m,b,c,h,e\rangle$ & $\langle p,o\rangle$ \\
\hline
$b$ & $p$ & $d$ & $m$ & $\langle b,a\rangle$ & $ \langle p,o,n,i,l,k,j,g,h,c,f,e\rangle$ & $\langle d,m\rangle$ \\
\hline
$b$ & $p$ & $f$ & $c$ & $\langle b,a\rangle$ & $ \langle p,m,d,o,n,i,l,k,j,g,h,e\rangle$ & $\langle f,c\rangle$ \\
\hline
$b$ & $p$ & $j$ & $i$ & $\langle b,a\rangle$ & $ \langle p,m,d,o,n,k,l,g,h,c,f,e\rangle$ & $\langle j,i\rangle$ \\
\hline
$b$ & $p$ & $n$ & $o$ & $\langle b,m,d,a\rangle$ & $ \langle p,i,l,k,j,g,h,c,f,e\rangle$ & $\langle n,o\rangle$ \\
\hline
$j$ & $h$ & $b$ & $c$ & $\langle j,k,l,i,n,o,p,m,d,a\rangle$ & $ \langle h,g,f,e\rangle$ & $\langle b,c\rangle$ \\
\hline
$j$ & $h$ & $f$ & $g$ & $\langle j,k,l,i,n,o,p,m,d,c,b,a\rangle$ & $ \langle h,e\rangle$ & $\langle f,g\rangle$ \\
\hline
$j$ & $h$ & $l$ & $k$ & $\langle j,i,n,o,p,m,d,c,b,a\rangle$ & $ \langle h,g,f,e\rangle$ & $\langle l,k\rangle$ \\
\hline
\end{tabular}
\end{center}

\begin{center}
\begin{tabular}{|l|l|l|l|l|l|l|}
\hline
$s_1$ & $s_2$ & $s_3$ & $t_3$ &$P_1[s_1,t_1]$ & $P_2[s_2,t_2]$ & $P_3[s_3,t_3]$ \\
\hline
$j$ & $h$ & $n$ & $o$ & $\langle j,k,l,i,p,m,d,c,b,a\rangle$ & $ \langle h,g,f,e\rangle$ & $\langle n,o\rangle$ \\
\hline
$f$ & $l$ & $b$ & $o$ & $\langle f,g,j,k,n,i,p,m,d,c,h,a\rangle$ & $ \langle l,e\rangle$ & $\langle b,o\rangle$ \\
\hline
$f$ & $l$ & $h$ & $g$ & $\langle f,c,d,o,p,m,b,a\rangle$ & $ \langle l,i,n,k,j,e\rangle$ & $\langle h,g\rangle$ \\
\hline
$f$ & $l$ & $j$ & $k$ & $\langle f,c,d,o,p,m,b,a\rangle$ & $ \langle l,g,h,e\rangle$ & $\langle j,i,n,k\rangle$ \\
\hline
$f$ & $l$ & $n$ & $i$ & $\langle f,c,d,o,p,m,b,a\rangle$ & $ \langle l,k,j,g,h,e\rangle$ & $\langle n,i\rangle$ \\
\hline
$n$ & $h$ & $b$ & $c$ & $\langle n,o,d,m,p,i,l,k,j,g,f,a\rangle$ & $ \langle h,e\rangle$ & $\langle b,c\rangle$ \\
\hline
$n$ & $h$ & $f$ & $c$ & $\langle n,i,p,o,d,m,b,a\rangle$ & $ \langle h,g,l,k,j,e\rangle$ & $\langle f,c\rangle$ \\
\hline
$n$ & $h$ & $j$ & $k$ & $\langle n,i,o,o,d,m,b,a\rangle$ & $ \langle h,c,f,e\rangle$ & $\langle j,g,l,k\rangle$ \\
\hline
$n$ & $h$ & $p$ & $k$ & $\langle n,o,d,m,b,a\rangle$ & $ \langle h,c,f,g,j,e\rangle$ & $\langle p,i,l,k\rangle$ \\
\hline
$f$ & $p$ & $b$ & $m$ & $\langle f,g,h,c,d,a\rangle$ & $ \langle p,o,n,i,l,k,j,e\rangle$ & $\langle b,m\rangle$ \\
\hline
$f$ & $p$ & $h$ & $g$ & $\langle f,c,b,m,d,a\rangle$ & $ \langle p,o,n,i,l,k,j,e\rangle$ & $\langle h,g\rangle$ \\
\hline
$f$ & $p$ & $j$ & $k$ & $\langle f,c,b,m,d,a\rangle$ & $ \langle p,o,n,i,l,g,h,e\rangle$ & $\langle j,k\rangle$ \\
\hline
$f$ & $p$ & $n$ & $o$ & $\langle f,c,b,m,d,a\rangle$ & $ \langle p,i,l,k.j,g,h,e\rangle$ & $\langle n,o\rangle$ \\
\hline
$n$ & $l$ & $b$ & $m$ & $\langle n,i,p,o,d,c,f,g,h,a\rangle$ & $ \langle l,k,j,e\rangle$ & $\langle b,m\rangle$ \\
\hline
$n$ & $l$ & $f$ & $g$ & $\langle n,k,j,i,p,o,d,m,b,c,h,a\rangle$ & $ \langle l,e\rangle$ & $\langle f,g\rangle$ \\
\hline
$n$ & $l$ & $j$ & $k$ & $\langle n,i,p,o,d,m,b,c,f,g,h,a\rangle$ & $ \langle l,e\rangle$ & $\langle j,k\rangle$ \\
\hline
$n$ & $l$ & $p$ & $i$ & $\langle n,o,d,m,b,c,f,g,h,a\rangle$ & $ \langle l,k,j,e\rangle$ & $\langle p,i\rangle$ \\
\hline
$j$ & $p$ & $b$ & $o$ & $\langle j,g,h,c,f,a\rangle$ & $ \langle p,i,n,k,l,e\rangle$ & $\langle b,m,d,o\rangle$ \\
\hline
$j$ & $p$ & $f$ & $g$ & $\langle j,k,l,i,n,o,d,a\rangle$ & $ \langle p,m,b,c,h,e\rangle$ & $\langle f,g\rangle$ \\
\hline
$j$ & $p$ & $l$ & $i$ & $\langle j,k,n,o,b,a\rangle$ & $ \langle p,m,d,c,h,g,f,e\rangle$ & $\langle l,i\rangle$ \\
\hline
$j$ & $p$ & $n$ & $k$ & $\langle j,g,h,a\rangle$ & $ \langle p,o,d,m,b,c,f,e\rangle$ & $\langle n,i,l,k\rangle$ \\
\hline

\end{tabular}
\end{center}

\medskip

\begin{center}
\texttt{Table 3. $t_1=(0,0)$ and $t_2=(0,2)$}
\medskip

\begin{tabular}{|l|l|l|l|l|l|l|}
\hline
$s_1$ & $s_2$ & $s_3$ & $t_3$ &$P_1[s_1,t_1]$ & $P_2[s_2,t_2]$ & $P_3[s_3,t_3]$ \\
\hline

$b$ & $d$ & $f$ & $e$ & $\langle b,c,h,a\rangle$ & $ \langle d,m,p,o,n,k,j,g,l,i\rangle$ & $\langle f,e\rangle$ \\
\hline
$b$ & $d$ & $j$ & $g$ & $\langle b,c,f,e,h,a\rangle$ & $ \langle d,m,p,o,n,k,l,i\rangle$ & $\langle j,g\rangle$ \\
\hline
$b$ & $d$ & $n$ & $o$ & $\langle b,c,f,e,h,a\rangle$ & $ \langle d,m,p,k,j,g,l,i\rangle$ & $\langle n,o\rangle$ \\
\hline
$f$ & $h$ & $b$ & $c$ & $\langle f,g,l,k,n,o,p,m,d,a\rangle$ & $ \langle h,e,j,i\rangle$ & $\langle b,c\rangle$ \\

$f$ & $h$ & $j$ & $e$ & $\langle f,c,b,m,n,k,p,o,d,a\rangle$ & $ \langle h,g,l,i\rangle$ & $\langle j,e\rangle$ \\
\hline
$f$ & $h$ & $n$ & $o$ & $\langle f,g,j,k,p,m,d,c,b,a\rangle$ & $ \langle h,e,l,i\rangle$ & $\langle n,o\rangle$ \\
\hline
$j$ & $l$ & $b$ & $c$ & $\langle j,g,h,e,f,a\rangle$ & $ \langle l,k,p,m,d,o,n,i\rangle$ & $\langle b,c\rangle$ \\
\hline
$j$ & $l$ & $f$ & $e$ & $\langle j,g,h,c,b,a\rangle$ & $ \langle l,k,p,m,d,o,n,i\rangle$ & $\langle f,e\rangle$ \\
\hline
$j$ & $l$ & $n$ & $m$ & $\langle j,g,f,e,h,c,b,a\rangle$ & $ \langle l,k,p,i\rangle$ & $\langle n,o,d,m\rangle$ \\
\hline
$n$ & $p$ & $b$ & $o$ & $\langle n,k,l,g,j,e,h,c,f,a\rangle$ & $ \langle p,i\rangle$ & $\langle b,m,d,o\rangle$ \\
\hline
$n$ & $p$ & $f$ & $e$ & $\langle n,o,b,a\rangle$ & $ \langle p,m,d,c,h,g,j,k,l,i\rangle$ & $\langle f,e\rangle$ \\
\hline
$n$ & $p$ & $j$ & $k$ & $\langle n,o,b,a\rangle$ & $ \langle p,m,d,c,h,e,f,g,l,i\rangle$ & $\langle j,k\rangle$ \\
\hline
$f$ & $d$ & $b$ & $c$ & $\langle f,e,h,a\rangle$ & $ \langle d,o,p,m,n,k,l.g,k,i\rangle$ & $\langle b,c\rangle$ \\
\hline
$f$ & $d$ & $h$ & $e$ & $\langle f,a\rangle$ & $ \langle d,c,b,m,p,o,n,k,j,g,l,i\rangle$ & $\langle h,e\rangle$ \\
\hline
$f$ & $d$ & $j$ & $g$ & $\langle f,e,h,c,b,a\rangle$ & $ \langle d,o,p,m,n,k,l,i\rangle$ & $\langle j,g\rangle$ \\
\hline
$f$ & $d$ & $n$ & $m$ & $\langle f,e,h,c,b,a\rangle$ & $ \langle d,o,p,k,l,g,j,i\rangle$ & $\langle n,m\rangle$ \\
\hline
$b$ & $h$ & $d$ & $c$ & $\langle b,m,n,o,p,k,l,e,f,a\rangle$ & $ \langle h,g,j,i\rangle$ & $\langle d,c\rangle$ \\
\hline
\end{tabular}
\end{center}

\begin{center}
\begin{tabular}{|l|l|l|l|l|l|l|}
\hline
$s_1$ & $s_2$ & $s_3$ & $t_3$ &$P_1[s_1,t_1]$ & $P_2[s_2,t_2]$ & $P_3[s_3,t_3]$ \\
\hline
$b$ & $h$ & $f$ & $c$ & $\langle b,o,p,m,d,a\rangle$ & $ \langle h,e,j,g,l,k,n,i\rangle$ & $\langle f,c\rangle$ \\
\hline
$b$ & $h$ & $j$ & $g$ & $\langle b,o,p,m,d,a\rangle$ & $ \langle h,c,f,e,l,k,n,i\rangle$ & $\langle j,g\rangle$ \\
\hline
$b$ & $h$ & $n$ & $k$ & $\langle b,m,p,o,d,a\rangle$ & $ \langle h,c,f,e,j,g,l,i\rangle$ & $\langle n,k\rangle$ \\
\hline
$j$ & $d$ & $b$ & $c$ & $\langle j,e,f,g,h,a\rangle$ & $ \langle d,m,p,o,n,k,l,i\rangle$ & $\langle b,c\rangle$ \\
\hline
$j$ & $d$ & $f$ & $e$ & $\langle j,g,h,c,b,a\rangle$ & $ \langle d,m,p,o,n,k,l,i\rangle$ & $\langle f,e\rangle$ \\
\hline
$j$ & $d$ & $l$ & $k$ & $\langle j,e,h,g,f,a\rangle$ & $ \langle d,c,b,m,n,o,p,i\rangle$ & $\langle l,k\rangle$ \\
\hline
$j$ & $d$ & $n$ & $o$ & $\langle j,g,f,e,h,c,b,a\rangle$ & $ \langle d,m,p,k,l,i\rangle$ & $\langle n,o\rangle$ \\
\hline
$b$ & $l$ & $d$ & $c$ & $\langle b,m,n,o,p,k,j,g,h,e,f,a\rangle$ & $ \langle l,i\rangle$ & $\langle d,c\rangle$ \\
\hline
$b$ & $l$ & $f$ & $e$ & $\langle b,m,n,o,p,k,j,g,h,c,d,a\rangle$ & $ \langle l,i\rangle$ & $\langle f,e\rangle$ \\
\hline
$b$ & $l$ & $j$ & $g$ & $\langle b,o,p,m,d,c,h,e,f,a\rangle$ & $ \langle l,n,k,i\rangle$ & $\langle j,g\rangle$ \\
\hline
$b$ & $l$ & $n$ & $o$ & $\langle b,c,d,m,p,k,j,g,h,e,f,a\rangle$ & $ \langle l,i\rangle$ & $\langle n,o\rangle$ \\
\hline
$n$ & $d$ & $b$ & $m$ & $\langle n,k,l,g,j,e,f,c,h,a\rangle$ & $ \langle d,o,p,i\rangle$ & $\langle b,m\rangle$ \\
\hline
$n$ & $d$ & $f$ & $e$ & $\langle n,m,p,o,b,a\rangle$ & $ \langle d,c,h,g,j,k,l,i\rangle$ & $\langle f,e\rangle$ \\
\hline
$n$ & $d$ & $j$ & $k$ & $\langle n,m,p,o,b,a\rangle$ & $ \langle d,c,f,g,h,e,l,i\rangle$ & $\langle j,k\rangle$ \\
\hline
$n$ & $d$ & $p$ & $m$ & $\langle n,o,b,a\rangle$ & $ \langle d,c,h,e,f,g,j,k,l,i\rangle$ & $\langle p,m\rangle$ \\
\hline
$b$ & $p$ & $d$ & $m$ & $\langle b,o,n,k,l,g,j,e,h,c,f,a\rangle$ & $ \langle p,i\rangle$ & $\langle d,m\rangle$ \\
\hline
$b$ & $p$ & $f$ & $c$ & $\langle b,o,n,m,d,a\rangle$ & $ \langle p,k,l,e,h,g,j,i\rangle$ & $\langle f,c\rangle$ \\
\hline
$b$ & $p$ & $j$ & $g$ & $\langle b,m,d,o,n,k,l,e,f,c,h,a\rangle$ & $ \langle p,i\rangle$ & $\langle j,g\rangle$ \\
\hline
$b$ & $p$ & $n$ & $o$ & $\langle b,m,d,c,h,e,l,k,j,g,f,a\rangle$ & $ \langle p,i\rangle$ & $\langle n,o\rangle$ \\
\hline
$j$ & $h$ & $b$ & $c$ & $\langle j,k,n,o,p,m,d,a\rangle$ & $ \langle h,e,f,g,l,i\rangle$ & $\langle b,c\rangle$ \\
\hline
$j$ & $h$ & $f$ & $e$ & $\langle j,k,n,o,p,m,d,c,b,a\rangle$ & $ \langle h,g,l,i\rangle$ & $\langle f,e\rangle$ \\
\hline
$j$ & $h$ & $l$ & $g$ & $\langle j,k,p,o,b,a\rangle$ & $ \langle h,e,f,c,d,m,n,i\rangle$ & $\langle l,g\rangle$ \\
\hline
$j$ & $h$ & $n$ & $o$ & $\langle j,k,p,m,d,c,b,a\rangle$ & $ \langle h,g,f,e,l,i\rangle$ & $\langle n,o\rangle$ \\
\hline
$f$ & $l$ & $b$ & $m$ & $\langle f,e,j,g,g,c,d,a\rangle$ & $ \langle l,k,p,o,n,i\rangle$ & $\langle b,m\rangle$ \\
\hline
$f$ & $l$ & $h$ & $e$ & $\langle f,c,d,m,b,a\rangle$ & $ \langle l,g,j,k,p,o,n,i\rangle$ & $\langle h,e\rangle$ \\
\hline
$f$ & $l$ & $j$ & $g$ & $\langle f,e,h,c,d,m,b,a\rangle$ & $ \langle l,k,p,o,n,i\rangle$ & $\langle j,g\rangle$ \\
\hline
$f$ & $l$ & $n$ & $o$ & $\langle f,g,j,e,h,c,d,m,b,a\rangle$ & $ \langle l,k,p,i\rangle$ & $\langle n,o\rangle$ \\
\hline
$n$ & $h$ & $b$ & $c$ & $\langle n,o,p,m,d,a\rangle$ & $ \langle h,e,f,g,j,k,l,i\rangle$ & $\langle b,c\rangle$ \\
\hline
$n$ & $h$ & $f$ & $e$ & $\langle n,o,p,m,d,c,b,a\rangle$ & $ \langle h,g,j,k,l,i\rangle$ & $\langle f,e\rangle$ \\
\hline
$n$ & $h$ & $j$ & $k$ & $\langle n,o,p,m,d,c,b,a\rangle$ & $ \langle h,g,f,e,l,i\rangle$ & $\langle j,k\rangle$ \\
\hline
$n$ & $h$ & $p$ & $o$ & $\langle n,m,d,c,b,a\rangle$ & $ \langle h,e,f,g,j,k,l,i\rangle$ & $\langle p,o\rangle$ \\
\hline
$f$ & $p$ & $b$ & $m$ & $\langle f,g,j,e,h,c,d,a\rangle$ & $ \langle p,o,n,k,l,i\rangle$ & $\langle b,m\rangle$ \\
\hline
$f$ & $p$ & $h$ & $e$ & $\langle f,c,d,m,b,a\rangle$ & $ \langle p,o,n,k,l,g,j,i\rangle$ & $\langle h,e\rangle$ \\
\hline
$f$ & $p$ & $j$ & $g$ & $\langle f,e,h,c,d,m,b,a\rangle$ & $ \langle p,o,n,k,l,i\rangle$ & $\langle j,g\rangle$ \\
\hline
$f$ & $p$ & $n$ & $o$ & $\langle f,e,j,g,h,c,d,m,b,a\rangle$ & $ \langle p,k,l,i\rangle$ & $\langle n,o\rangle$ \\
\hline
$n$ & $l$ & $b$ & $m$ & $\langle n,k,p,o,d,c,h,a\rangle$ & $ \langle l,e,f,g,j,i\rangle$ & $\langle b,m\rangle$ \\
\hline
$n$ & $l$ & $f$ & $c$ & $\langle n,k,p,o,d,m,b,a\rangle$ & $ \langle l,e,h,g,j,i\rangle$ & $\langle f,c\rangle$ \\
\hline
$n$ & $l$ & $j$ & $g$ & $\langle n,k,p,o,d,m,b,c,f,e,h,a\rangle$ & $ \langle l,i\rangle$ & $\langle j,g\rangle$ \\
\hline
$n$ & $l$ & $p$ & $k$ & $\langle n,o,d,m,b,c,f,e,h,a\rangle$ & $ \langle l,g,j,i\rangle$ & $\langle p,k\rangle$ \\
\hline
$j$ & $p$ & $b$ & $m$ & $\langle j,k,l,g,f,e,h,c,d,a\rangle$ & $ \langle p,o,n,i\rangle$ & $\langle b,m\rangle$ \\
\hline
$j$ & $p$ & $f$ & $e$ & $\langle j,k,l,g,h,c,d,m,b,a\rangle$ & $ \langle p,o,n,i\rangle$ & $\langle f,e\rangle$ \\
\hline
$j$ & $p$ & $l$ & $k$ & $\langle j,g,f,e,h,c,d,m,b,a\rangle$ & $ \langle p,o,n,i\rangle$ & $\langle l,k\rangle$ \\
\hline
$j$ & $p$ & $n$ & $o$ & $\langle j,k,l,g,f,e,h,c,b,m,d,a\rangle$ & $ \langle p,i\rangle$ & $\langle n,o\rangle$ \\
\hline
\end{tabular}
\end{center}

\newpage

\begin{center}
\texttt{Table 4. $t_1=(0,0)$ and $t_2=(0,3)$}
\medskip

\begin{tabular}{|l|l|l|l|l|l|l|}
\hline
$s_1$ & $s_2$ & $s_3$ & $t_3$ &$P_1[s_1,t_1]$ & $P_2[s_2,t_2]$ & $P_3[s_3,t_3]$ \\
\hline

$b$ & $d$ & $f$ & $e$ & $\langle b,a\rangle$ & $ \langle d,c,h,g,j,k,l,i,n,o,p,m\rangle$ & $\langle f,e\rangle$ \\
\hline
$b$ & $d$ & $j$ & $k$ & $\langle b,a\rangle$ & $ \langle d,c,h,e,f,g,l,i,n,o,p,m\rangle$ & $\langle j,k\rangle$ \\
\hline
$b$ & $d$ & $n$ & $o$ & $\langle b,a\rangle$ & $ \langle d,c,h,e,f,g,j,i,l,k,p,m\rangle$ & $\langle n,o\rangle$ \\
\hline
$f$ & $h$ & $b$ & $o$ & $\langle f,c,d,a\rangle$ & $ \langle h,e,j,g,l,k,n,i,p,m\rangle$ & $\langle b,o\rangle$ \\
\hline
$f$ & $h$ & $j$ & $e$ & $\langle f,c,b,a\rangle$ & $ \langle h,g,l,k,n,i,p,o,d,m\rangle$ & $\langle j,e\rangle$ \\
\hline
$f$ & $h$ & $n$ & $k$ & $\langle f,c,b,a\rangle$ & $ \langle h,e,j,g,l,i,p,o,d,m\rangle$ & $\langle n,k\rangle$ \\
\hline
$j$ & $l$ & $b$ & $o$ & $\langle j,g,f,e,h,c,d,a\rangle$ & $ \langle l,k,n,i,p,m\rangle$ & $\langle b,o\rangle$ \\
\hline
$j$ & $l$ & $f$ & $g$ & $\langle j,e,h,c,d,a\rangle$ & $ \langle l,k,n,i,p,o,b,m\rangle$ & $\langle f,g\rangle$ \\
\hline
$j$ & $l$ & $n$ & $i$ & $\langle j,g,f,e,h,c,d,a\rangle$ & $ \langle l,k,p,o,b,m\rangle$ & $\langle n,i\rangle$ \\
\hline
$n$ & $p$ & $b$ & $o$ & $\langle n,i,l,k,j,g,f,e,h,c,d,a\rangle$ & $ \langle p,m\rangle$ & $\langle b,o\rangle$ \\
\hline
$n$ & $p$ & $f$ & $g$ & $\langle n,i,l,k,j,e,h,c,d,a\rangle$ & $ \langle p,o,b,m\rangle$ & $\langle f,g\rangle$ \\
\hline
$n$ & $p$ & $j$ & $k$ & $\langle n,i,l,g,f,e,h,c,d,a\rangle$ & $ \langle p,o,b,m\rangle$ & $\langle j,k\rangle$ \\
\hline
$f$ & $d$ & $b$ & $c$ & $\langle f,g,j,e,h,a\rangle$ & $ \langle d,o,n,k,l,i,p,m\rangle$ & $\langle b,c\rangle$ \\
\hline
$f$ & $d$ & $h$ & $e$ & $\langle f,c,b,a\rangle$ & $ \langle d,o,n,k,j,g,l,i,p,m\rangle$ & $\langle h,e\rangle$ \\
\hline
$f$ & $d$ & $j$ & $g$ & $\langle f,e,h,c,b,a\rangle$ & $ \langle d,o,n,k,l,i,p,m\rangle$ & $\langle j,g\rangle$ \\
\hline
$f$ & $d$ & $n$ & $i$ & $\langle f,g,l,k,j,e,h,c,b,a\rangle$ & $ \langle d,o,p,m\rangle$ & $\langle n,i\rangle$ \\
\hline
$b$ & $h$ & $d$ & $c$ & $\langle b,a\rangle$ & $ \langle h,e,f,g,j,k,l,i,n,o,p,m\rangle$ & $\langle d,c\rangle$ \\
\hline
$b$ & $h$ & $f$ & $e$ & $\langle b,c,d,a\rangle$ & $ \langle h,g,j,k,l,i,n,o,p,m\rangle$ & $\langle f,e\rangle$ \\
\hline
$b$ & $h$ & $j$ & $i$ & $\langle b,c,d,a\rangle$ & $ \langle h,e,f,g,l,i,n,o,p,m\rangle$ & $\langle j,i\rangle$ \\
\hline
$b$ & $h$ & $n$ & $o$ & $\langle b,c,d,a\rangle$ & $ \langle h,e,f,g,j,i,l,k,p,m\rangle$ & $\langle n,o\rangle$ \\
\hline
$b$ & $l$ & $d$ & $o$ & $\langle b,c,f,g,j,e,h,a\rangle$ & $ \langle l,k,n,i,p,m\rangle$ & $\langle d,o\rangle$ \\
\hline
$b$ & $l$ & $f$ & $e$ & $\langle b,o,d,c,h,a\rangle$ & $ \langle l,g,j,k,n,i,p,m\rangle$ & $\langle f,e\rangle$ \\
\hline
$b$ & $l$ & $j$ & $g$ & $\langle b,o,d,c,f,e,h,a\rangle$ & $ \langle l,k,n,i,p,m\rangle$ & $\langle j,g\rangle$ \\
\hline
$b$ & $l$ & $n$ & $o$ & $\langle b,c,f,g,j,e,h,a\rangle$ & $ \langle l,k,j,i,p,m\rangle$ & $\langle n,o\rangle$ \\
\hline
$j$ & $d$ & $b$ & $c$ & $\langle j,k,l,g,f,e,h,a\rangle$ & $ \langle d,o,n,i,p,m\rangle$ & $\langle b,c\rangle$ \\
\hline
$j$ & $d$ & $f$ & $g$ & $\langle j,k,l,e,h,c,b,a\rangle$ & $ \langle d,o,n,i,p,m\rangle$ & $\langle f,g\rangle$ \\
\hline
$j$ & $d$ & $l$ & $k$ & $\langle j,g,f,e,h,c,b,a\rangle$ & $ \langle d,o,n,i,p,m\rangle$ & $\langle l,k\rangle$ \\
\hline
$j$ & $d$ & $n$ & $k$ & $\langle j,g,l,e,f,c,h,a\rangle$ & $ \langle d,o,b,m\rangle$ & $\langle n,i,p,k\rangle$ \\
\hline
$n$ & $d$ & $b$ & $o$ & $\langle n,i,p,k,l,g,j,e,f,c,h,a\rangle$ & $ \langle d,m\rangle$ & $\langle b,o\rangle$ \\
\hline
$n$ & $d$ & $f$ & $g$ & $\langle n,o,p,i,l,k,j,e,h,c,b,a\rangle$ & $ \langle d,m\rangle$ & $\langle f,g\rangle$ \\
\hline
$n$ & $d$ & $j$ & $k$ & $\langle n,o,p,i,l,g,f,e,h,c,b,a\rangle$ & $ \langle d,m\rangle$ & $\langle j,k\rangle$ \\
\hline
$n$ & $d$ & $p$ & $o$ & $\langle n,i,l,k,j,g,f,e,h,c,b,a\rangle$ & $ \langle d,m\rangle$ & $\langle p,o\rangle$ \\
\hline
$b$ & $p$ & $d$ & $o$ & $\langle b,c,f,g,j,e,h,a\rangle$ & $ \langle p,k,l,i,n,m\rangle$ & $\langle d,o\rangle$ \\
\hline
$b$ & $p$ & $f$ & $e$ & $\langle b,o,d,c,h,a\rangle$ & $ \langle p,k,j,g,l,i,n,m\rangle$ & $\langle f,e\rangle$ \\
\hline
$b$ & $p$ & $j$ & $g$ & $\langle b,o,d,c,f,e,h,a\rangle$ & $ \langle p,k,l,i,n,m\rangle$ & $\langle j,g\rangle$ \\
\hline
$b$ & $p$ & $n$ & $k$ & $\langle b,o,d,c,f,g,j,e,h,a\rangle$ & $ \langle p,m\rangle$ & $\langle n,i,l,k\rangle$ \\
\hline
$f$ & $l$ & $b$ & $o$ & $\langle f,g,j,e,h,c,d,a\rangle$ & $ \langle l,k,n,i,p,m\rangle$ & $\langle b,o\rangle$ \\
\hline
$f$ & $l$ & $h$ & $e$ & $\langle f,c,d,o,b,a\rangle$ & $ \langle l,g,j,k,n,i,p,m\rangle$ & $\langle h,e\rangle$ \\
\hline
\end{tabular}
\end{center}

\begin{center}
\begin{tabular}{|l|l|l|l|l|l|l|}
\hline
$s_1$ & $s_2$ & $s_3$ & $t_3$ &$P_1[s_1,t_1]$ & $P_2[s_2,t_2]$ & $P_3[s_3,t_3]$ \\
\hline
$f$ & $l$ & $j$ & $g$ & $\langle f,e,h,c,d,o,b,a\rangle$ & $ \langle l,k,n,i,p,m\rangle$ & $\langle j,g\rangle$ \\
\hline
$f$ & $l$ & $n$ & $i$ & $\langle f,g,j,e,h,c,d,o,b,a\rangle$ & $ \langle l,k,p,m\rangle$ & $\langle n,i\rangle$ \\
\hline
$j$ & $h$ & $b$ & $o$ & $\langle j,e,f,c,d,a\rangle$ & $ \langle h,g,l,k,n,i,p,m\rangle$ & $\langle b,o\rangle$ \\
\hline
$j$ & $h$ & $f$ & $e$ & $\langle j,i,n,o,b,c,d,a\rangle$ & $ \langle h,g,l,k,p,m\rangle$ & $\langle f,e\rangle$ \\
\hline
$j$ & $h$ & $l$ & $g$ & $\langle j,k,n,i,p,o,d,a\rangle$ & $ \langle h,e,f,c,b,m\rangle$ & $\langle l,g\rangle$ \\
\hline
$j$ & $h$ & $n$ & $i$ & $\langle j,e,f,c,d,a\rangle$ & $ \langle h,g,l,k,p,o,b,m\rangle$ & $\langle n,i\rangle$ \\
\hline
$f$ & $p$ & $b$ & $o$ & $\langle f,g,l,i,n,k,j,e,h,c,d,a\rangle$ & $ \langle p,m\rangle$ & $\langle b,o\rangle$ \\
\hline
$f$ & $p$ & $h$ & $e$ & $\langle f,c,b,a\rangle$ & $ \langle p,k,j,g,l,i,n,o,d,m\rangle$ & $\langle h,e\rangle$ \\
\hline
$f$ & $p$ & $j$ & $g$ & $\langle f,e,h,c,b,o,d,a\rangle$ & $ \langle p,k,l,i,n,m\rangle$ & $\langle j,g\rangle$ \\
\hline
$f$ & $p$ & $n$ & $i$ & $\langle f,g,l,k,j,e,h,c,b,o,d,a\rangle$ & $ \langle p,m\rangle$ & $\langle n,i\rangle$ \\
\hline
$n$ & $h$ & $b$ & $o$ & $\langle n,i,p,k,l,e,j,g,f,a\rangle$ & $ \langle h,c,d,m\rangle$ & $\langle b,o\rangle$ \\
\hline
$n$ & $h$ & $f$ & $e$ & $\langle n,o,d,c,b,a\rangle$ & $ \langle h,g,j,i,l,k,p,m\rangle$ & $\langle f,e\rangle$ \\
\hline
$n$ & $h$ & $j$ & $e$ & $\langle n,i,p,k,l,g,f,a\rangle$ & $ \langle h,c,b,o,d,m\rangle$ & $\langle j,e\rangle$ \\
\hline
$n$ & $h$ & $p$ & $i$ & $\langle n,k,l,e,j,g,f,a\rangle$ & $ \langle h,c,d,o,b,m\rangle$ & $\langle p,i\rangle$ \\
\hline
$j$ & $p$ & $b$ & $o$ & $\langle j,g,f,e,h,c,d,a\rangle$ & $ \langle p,i,l,k,n,m\rangle$ & $\langle b,o\rangle$ \\
\hline
$j$ & $p$ & $f$ & $g$ & $\langle j,e,h,c,d,a\rangle$ & $ \langle p,k,l,i,n,o,b,m\rangle$ & $\langle f,g\rangle$ \\
\hline
$j$ & $p$ & $l$ & $k$ & $\langle j,g,f,e,h,c,d,a\rangle$ & $ \langle p,i,n,o,b,m\rangle$ & $\langle l,k\rangle$ \\
\hline
$j$ & $p$ & $n$ & $i$ & $\langle j,k,l,g,f,e,h,c,d,a\rangle$ & $ \langle p,o,b,m\rangle$ & $\langle n,i\rangle$ \\
\hline
$n$ & $l$ & $b$ & $o$ & $\langle n,k,j,g,f,e,h,c,d,a\rangle$ & $ \langle l,i,p,m\rangle$ & $\langle b,o\rangle$ \\
\hline
$n$ & $l$ & $f$ & $g$ & $\langle n,k,j,e,h,c,d,a\rangle$ & $ \langle l,i,p,o,b,m\rangle$ & $\langle f,g\rangle$ \\
\hline
$n$ & $l$ & $j$ & $k$ & $\langle n,i,p,o,d,a\rangle$ & $ \langle l,g,f,e,h,c,b,m\rangle$ & $\langle j,k\rangle$ \\
\hline
$n$ & $l$ & $p$ & $i$ & $\langle n,o,d,a\rangle$ & $ \langle l,k,j,g,f,e,h,c,b,m\rangle$ & $\langle p,i\rangle$ \\
\hline
\end{tabular}
\end{center}
}

\end{document}